%% file: snake.tex
\documentclass[11pt]{article}

\usepackage[a4paper]{geometry}
\usepackage[LSF, X2, T1]{fontenc}
\usepackage{indentfirst}
\usepackage{hyperref}
\usepackage[anythingbreaks]{breakurl}
\usepackage{amssymb, amsmath}
\usepackage{graphicx}
\usepackage[labelsep=none]{caption}
\usepackage[labelformat=simple]{subcaption}

\hypersetup{pdfstartview=FitH, pdfpagemode=UseNone}

\allowdisplaybreaks

\newcounter{question}
\newenvironment{question}{\refstepcounter{question}\textbf{Question \thequestion.}\begin{em}}{\end{em}}
\newenvironment{question*}[1]{\textbf{Question #1.}\begin{em}}{\end{em}}
\newcounter{lemma}
\newenvironment{lemma}{\refstepcounter{lemma}\textbf{Lemma \thelemma.}\begin{em}}{\end{em}}
\newcounter{proposition}
\newenvironment{proposition}{\refstepcounter{proposition}\textbf{Proposition \theproposition.}\begin{em}}{\end{em}}
\newcounter{theorem}
\newenvironment{theorem}{\refstepcounter{theorem}\textbf{Theorem \thetheorem.}\begin{em}}{\end{em}}
\newenvironment{proof}{\emph{Proof.}}{$\square$}
\newenvironment{proof*}[1]{\emph{#1.}}{$\square$}

\def\spring{{\sim}}

\def\nosovka{\textrm{\usefont{X2}{cmr}{m}{n} \CYRBYUS}}

\makeatletter
\def\myking{\textrm{\usefont{LSF}{skaknew}{m}{n} \cfss@symking}}
\def\myknight{\textrm{\usefont{LSF}{skaknew}{m}{n} \cfss@symknight}}
\makeatother

\newcommand\kinggraph[1]{\mathcal{G}(\myking, #1)}
\newcommand\knightgraph[1]{\mathcal{G}(\myknight, #1)}
\newcommand\gridgraph[1]{\mathcal{G}(\square, #1)}

\begin{document}

\title{\textbf{Snake Paths in King and Knight Graphs}}
\author{Nikolai Beluhov}
\date{}

\maketitle

\begin{center} \parbox{352pt}{\setlength{\parindent}{12pt} \footnotesize \emph{Abstract}. A snake path in a graph $G$ is a path in $G$ which is also an induced subgraph of $G$. For all $n$, we find the greatest length of a snake path in the $n \times n$ king graph and we give a complete description of the paths which attain this greatest length. The even and odd cases behave very differently. We also estimate the greatest length of a snake path or cycle in the $m \times n$ knight graph, for all $m$ and $n$.} \end{center}

\input{snake-01-intro}

\input{snake-02-prelim}

\input{snake-03-king-even}

\input{snake-04-king-odd-i}

\input{snake-05-king-odd-ii}

\input{snake-06-knight}

\input{snake-07-further-i}

\input{snake-08-further-ii}

\input{snake-09-ack}

\input{snake-10-refs}
\input{snake-11-add}

\end{document}

%% file: snake-01-intro.tex
\section{Introduction} \label{intro}

Let $G$ be a simple graph. A \emph{snake path} in $G$ is a path in $G$ which is also an induced subgraph of $G$. Equivalently, a path $P$ in $G$ is a snake path when, for all vertices $u$ and $v$ of $P$, we have that $uv$ is an edge of $G$ if and only if it is an edge of $P$. Intuitively, a snake path never comes into contact with itself.

Snake paths are also known as \emph{induced paths} and \emph{chordless paths}.

A \emph{snake cycle} is defined similarly. Just like paths, snake cycles are alternatively called \emph{induced cycles} and \emph{chordless cycles}. Our focus will be mostly on snake paths, though we will touch upon snake cycles, too.

Given a graph $G$, some of the most natural questions we can ask about its snake paths are as follows:

\medskip

\begin{question*}{A} What is the greatest length of a snake path in $G$? \end{question*}

\medskip

Note that we measure the length of a path by the number of edges that it traverses, rather than the number of vertices that it visits.

The answer of Question \textbf{A} coincides with the greatest diameter of an induced subgraph of $G$.

The greatest length of a snake path in $G$ is also known as the \emph{induced detour number} of $G$, and the greatest length of a snake cycle in $G$ as the \emph{induced circumference} of $G$.

In the special case when $G$ is a hypercube graph, Question \textbf{A} is known as the \emph{snake-in-the-box problem}, and its analogue for cycles as the \emph{coil-in-the-box problem}. Both of these problems have been studied extensively.

\medskip

\begin{question*}{B} What is the structure of the longest snake paths in $G$? \end{question*}

\medskip

\begin{question*}{C} How many longest snake paths are there in $G$? \end{question*}

\medskip

Note that we formalise paths as subgraphs, rather than as sequences of vertices. Specifically, to us a ``path'' is a tree subgraph where all vertices are of degree at most two. The distinction between the two formalisations does not matter in most situations, but it does matter for enumeration. For example, to us $abc$ and $cba$ are the same path.

We study Questions \textbf{A}--\textbf{C} for certain graphs $G$ associated with chess pieces.

Given a chess piece $F$ and a board $A$, it is natural to consider the graph whose vertices are the cells of $A$ and whose edges correspond to all possible moves of $F$ on $A$. We proceed to formalise this notion for the king and the knight. Other chess pieces can be handled similarly.

To us, a \emph{cell} is an ordered pair of integers. Or, equivalently, an integer point in the plane.

Let $m$ and $n$ be positive integers. A \emph{board} $A$ of size $m \times n$, with $m$ rows and $n$ columns, is a set of cells of the form $I \times J$, where $I$ and $J$ are integer intervals with $|I| = n$ and $|J| = m$. The \emph{standard board} of size $m \times n$ has $I = [0; n - 1]$ and $J = [0; m - 1]$. Since all boards of the same size are translation copies of one another, sometimes we refer to ``the'' board of a certain size, meaning the standard board of that size.

Given a set of cells $S$, we define the \emph{king graph} on $S$, denoted $\kinggraph{S}$, to be the graph on vertex set $S$ where two distinct cells $a' = (x', y')$ and $a'' = (x'', y'')$ are joined by an edge if and only if $|x' - x''| \le 1$ and $|y' - y''| \le 1$. Since all king graphs on boards of the same size are isomorphic, sometimes we refer to ``the'' king graph of a certain size, meaning the king graph on the standard board of that size. For convenience, we also use the notation $\kinggraph{m \times n}$ for the king graph of size $m \times n$. Of course, $\kinggraph{m \times n}$ can also be viewed as the strong product of two paths with $m$ and $n$ vertices, respectively.

The \emph{knight graph} on $S$, denoted $\knightgraph{S}$, is defined similarly, except that the adjacency condition becomes $\{|x' - x''|, |y' - y''|\} = \{1, 2\}$ instead.

Dawson, in problem 187 of \cite{D}, considers $\knightgraph{8 \times 8}$ and presents a snake path of length $31$ as well as a snake cycle of length $32$.

Knuth, in exercise 172 of \cite{K}, discusses the longest snake paths and cycles of various chess piece graphs in the context of algorithmic generation. In particular, he determines that there are $16$ essentially distinct snake paths of the greatest length $31$ and $6$ essentially distinct snake cycles of the greatest length $31$ in $\kinggraph{8 \times 8}$; as well as an essentially unique snake path of the greatest length $33$ and $4$ essentially distinct snake cycles of the greatest length $32$ in $\knightgraph{8 \times 8}$. (Thus Dawson's path was not optimal, but his cycle was.)

Our main results are as follows:

\medskip

\begin{theorem} \label{king-even-path} Let $n$ be an even positive integer. Then the greatest length of a snake path in the king graph of size $n \times n$ is $n^2/2 - 1$. Furthermore, when $n \ge 6$, there are exactly $16n$ snake paths which attain this greatest length. \end{theorem}

\medskip

This count includes rotations and reflections. In Section \ref{king-even}, we will see that with $n \ge 6$ the number of essentially distinct longest snake paths is $2n$ when $n/2$ is even and $2n + 1$ when $n/2$ is odd.

Over the course of the proof of Theorem \ref{king-even-path}, we will also give a complete description of these paths. Roughly speaking, each one of them is shaped like a spiral.

That the total number of paths is given by such a nice formula is most likely only a happy coincidence. Because of the overall structure of the proof, we should expect to see a linear function of $n$ within each parity of $n/2$; however, there is no obvious reason \emph{a priori} to expect these two functions to coincide, or their constant terms to vanish.

\medskip

\begin{theorem} \label{king-odd-path} Let $n$ be an odd positive integer. Then the greatest length of a snake path in the king graph of size $n \times n$ is $(n^2 - 1)/2$. \end{theorem}

\medskip

Remarkably, the odd and even cases behave very differently. Despite the surface similarity between the upper bounds of Theorems \ref{king-even-path} and \ref{king-odd-path}, the former bound is straightforward while the latter one poses considerable difficulties.

In Section \ref{king-odd-ii}, we will add to Theorem \ref{king-odd-path} a complete description of the paths which attain the greatest length, stated in Theorem \ref{king-odd-structure}. However, the description is somewhat complicated, and relies on a long series of preceding definitions. Thus we do not reproduce it in the introduction. The gist is that each stamp-folding permutation of $\lceil n/2 \rceil$ elements yields two families of longest snake paths which share the same overall shape but differ from one another by some tiny aberrations.

Theorem \ref{king-odd-structure} does not imply an exact answer to Question \textbf{C}. In fact, because of the connection to the stamp-folding problem, it seems unlikely that such an answer would be feasible. Still, the theorem does yield some loose bounds. In particular, we will see (Proposition \ref{log-bound}) that the logarithm of the number of longest snake paths grows as $\Theta(n)$, in stark contrast to the even case.

Theorems \ref{king-even-path}, \ref{king-odd-path}, and \ref{king-odd-structure} together completely resolve the questions of the greatest length of a snake path and the structure of the longest snake paths in king graphs on square boards.

\medskip

\begin{theorem} \label{knight-all} Let $m$ and $n$ be positive integers. Then both the longest snake path and the longest snake cycle in the knight graph of size $m \times n$ are of length $mn/2 + \mathcal{O}(m + n)$. \end{theorem}

\medskip

Note that we do not specify the sign of the error term: Since $\mathcal{O}$-notation only bounds the absolute value of a function, the classes $mn/2 + \mathcal{O}(m + n)$ and $mn/2 - \mathcal{O}(m + n)$ consist of the same functions of $m$ and $n$. Same goes for the estimates in Section \ref{further-ii}.

Compared to the treatment of Question \textbf{A} in Theorems \ref{king-even-path} and \ref{king-odd-path}, with Theorem \ref{knight-all} we do not attempt to obtain an exact answer, and are content instead with an asymptotic estimate. On the bright side, this asymptotic estimate applies to cycles as well as paths, and it is valid on all rectangular boards.

It is worth noting that one step in our proof of Theorem \ref{knight-all} relies on computer help, and likely cannot be verified manually by a human mathematician. (Our proofs of Theorems \ref{king-even-path}, \ref{king-odd-path}, and \ref{king-odd-structure} are all human-friendly, though.)

The author obtained Theorems \ref{king-even-path}--\ref{near-hamiltonian} in 2018 after being introduced to the subject by Knuth, in connection with the aforementioned exercise 172. Subsequently, Theorems \ref{king-even-path}, \ref{knight-all}, and \ref{king-0mod4} were cited in a remark following the exercise's solution. (Strictly speaking, at that time the author derived Theorem \ref{king-even-path} in a form referring to the number of essentially distinct paths rather than the total number of paths. This is also how it was stated in~\cite{K}.)

Then, in 2020, the author proposed Theorem \ref{king-odd-path}, appropriately rephrased, as a mathematical olympiad problem for the Cyberspace Mathematical Competition. It was featured as problem 4 on day 1 of the contest. (The CMC was a one-off event intended to approximate the International Mathematical Olympiad. Perhaps the problem's difficulty was not an ideal match for the contest's format; out of 553 participants from 75 countries, only two made substantial progress on it.)

Theorems \ref{king-even-path}--\ref{king-odd-structure} offer an interesting illustration of how the nature of Questions \textbf{A}--\textbf{C} can change when we vary the underlying graph. In the setting of Theorem \ref{king-even-path}, Question \textbf{A} is straightforward while Questions \textbf{B} and \textbf{C} are manageable. For Theorems \ref{king-odd-path} and \ref{king-odd-structure}, both Questions \textbf{A} and \textbf{B} become significantly more complicated, while a closed-form answer to Question \textbf{C} is likely out of reach. Finally, in the setting of Theorem \ref{knight-all}, already with Question \textbf{A} it seems that an exact answer would be unfeasible, and even for our asymptotic estimate we find ourselves in need of machine help.

%% file: snake-02-prelim.tex
\section{Preliminaries} \label{prelim}

Before we continue, let us briefly list some useful notations and observations.

Given two cells $a' = (x', y')$ and $a'' = (x'', y'')$, we write $a' + a''$ for the cell $(x' + x'', y' + y'')$. Given a cell $a$ and a set of cells $S$, we write $a + S$ for the set of cells $\{a + b \mid b \in S\}$.

A \emph{symmetry} of a board $A$ is the restriction to $A$ of an isometry of the plane which preserves $A$. Two objects defined with reference to $A$, such as two sets of cells on $A$ or two graphs on $A$, are \emph{essentially distinct} (relative to $A$) when they are distinct under the symmetries of $A$.

The \emph{grid graph} on a set of cells $S$, denoted $\gridgraph{S}$, is defined similarly to the king and knight graphs on $S$, except that the adjacency condition becomes $\{|x' - x''|, |y' - y''|\} = \{0, 1\}$ instead. Or, equivalently, cells $a'$ and $a''$ must be at unit Euclidean distance from one another. Of course, $\gridgraph{m \times n}$ can also be viewed as the Cartesian product of two paths with $m$ and $n$ vertices, respectively.

For the edge joining two vertices $a$ and $b$ in a graph $G$, we write either $ab$ or $a$---$b$, whichever one reads better in the situation at hand.

We introduce shorthand notation for certain paths within grid and king graphs. Given two cells $a$ and $b$ of a grid or king graph $G$ such that $a$ and $b$ are in the same row or column, we write $a \spring b$ for the path in $G$ connecting $a$ and $b$ whose remaining cells are the ones between $a$ and $b$ in the corresponding row or column. For example, if $a = (x', y)$, $b = (x'', y)$, and $x' \le x''$, then $a \spring b = (x', y)$---$(x' + 1, y)$---$(x' + 2, y)$---$\cdots$---$(x'', y)$.

Let $P$ be a path in some graph $G$. Given a vertex $a$ of $G$, we write $a \in P$ for ``$P$ visits $a$''; given an edge $e$ of $G$, we write $e \in P$ for ``$P$ traverses $e$''; and, given a path $Q$ in $G$, we write $Q \subseteq P$ for ``$Q$ is a subpath of $P$''. We will not use these abbreviations too often, but the proofs of Lemmas \ref{spiral} and \ref{nice-ii} would be cumbersome to state without them.

Suppose, now, that $P$ is a snake path in $G$. When $ab'b''$ is a three-cycle in $G$, we have that $a \in P$ implies $b'b'' \not \in P$. Similarly, when $b$, $c'$, and $c''$ are three distinct neighbours of $a$ in $G$, we have that $c'ac'' \subseteq P$ implies $b \not \in P$. We will use these simple observations repeatedly throughout the paper.

A subset $S$ of the vertices of $G$ is \emph{$k$-independent} if, in the induced subgraph of $G$ on vertex set $S$, every vertex is of degree at most $k$. When $k = 2$, the induced subgraph itself is a \emph{pseudosnake} of $G$. Clearly, the number of vertices in the longest snake path or cycle of $G$ cannot exceed the number of vertices in its largest pseudosnake. When $G$ is finite, we define its \emph{pseudosnake density} to be the ratio of the greatest number of vertices in a pseudosnake of $G$ to the total number of vertices in $G$.

%% file: snake-03-king-even.tex
\section{King Graphs on Even Boards} \label{king-even}

In this section, we prove Theorem \ref{king-even-path}. Let $n$ be an even positive integer with $n = 2k$, let $A$ be the standard board of size $n \times n$, and let $G$ be the king graph on $A$.

For completeness, with $n = 2$ there are $6$ longest snake paths, of which $2$ are essentially distinct; and with $n = 4$ there are $28$ longest snake paths, of which $4$ are essentially distinct. From now on, let $n \ge 6$.

The argument we give for the upper bound of Theorem \ref{king-even-path} is not new; the special case $n = 8$ is in \cite{K}, and the general case does not pose any additional difficulties.

\medskip

\begin{proof*}{Proof of the optimisation part of Theorem \ref{king-even-path}} For the upper bound, partition $A$ into $k^2$ subboards of size $2 \times 2$ each. Since a snake path in $G$ can visit at most two cells within each subboard, its length cannot exceed $2k^2 - 1$.

\begin{figure}[ht] \centering \includegraphics{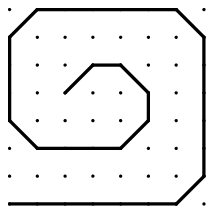} \caption{} \label{king-8x8} \end{figure}

For the lower bound, let $P_i$ be the path $(i, i) \spring (n - i - 2, i)$---$(n - i - 1, i + 1) \spring (n - i - 1, n - i - 2)$---$(n - i - 2, n - i - 1) \spring (i + 1, n - i - 1)$---$(i, n - i - 2) \spring (i, i + 3)$---$(i + 1, i + 2)$---$(i + 2, i + 2)$ for all even $i$ with $0 \le i \le k - 3$. When $k$ is even and $i = k - 2$, we define $P_{k - 2}$ in the same way, except that we stop at cell $(i, n - i - 2) = (k - 2, k)$; and, when $k$ is odd and $i = k - 1$, we also define $P_{k - 1}$ to be the path $(k - 1, k - 1)$---$(k, k)$. Then the concatenation of these paths is a snake path in $G$ of length $n^2/2 - 1$. \end{proof*}

\medskip

For example, Figure \ref{king-8x8} shows $n = 8$.

The enumeration part of Theorem \ref{king-even-path} will be somewhat more complicated.

Let $B$ be the standard board of size $k \times k$ and let $H$ be the grid graph on $B$.

For each cell $b$ of $B$, let $\Phi(b)$ denote the subboard $2b + [0; 1]^2$ of $A$, of size $2 \times 2$. (For convenience, if $b = (x, y)$, we also write simply $\Phi(x, y)$.) These subboards, taken over all cells $b$ of $B$, form a partitioning of $A$.

Let $P$ be a longest snake path in $G$. The proof of the optimisation part of Theorem~\ref{king-even-path} shows that $P$ visits exactly two cells within each subboard of $A$ of the form $\Phi(b)$.

Suppose that an edge of $P$ joins one cell of $\Phi(b')$ and one cell of $\Phi(b'')$, with $b' \neq b''$. We claim that $b'$ and $b''$ are then neighbours in $H$.

\begin{figure}[t] \centering \includegraphics{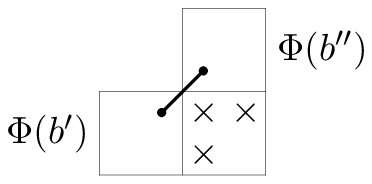} \caption{} \label{non-dia} \end{figure}

Indeed, if not, then $b'$ and $b''$ must be diagonally adjacent in $H$, without loss of generality with $b' + (1, 1) = b''$. So $P$ contains a subpath of the form $c'$---$(2b' + (1, 1))$---$2b''$---$c''$, where $c'$ is in $\Phi(b')$ and $c''$ is in $\Phi(b'')$. Since $P$ is a snake path, it follows that $P$ cannot visit any cells in the set $2(b' + (1, 0)) + \{(0, 0), (0, 1), (1, 1)\}$. Consequently, $P$ visits at most one cell of $\Phi(b' + (1, 0))$, a contradiction. (Figure \ref{non-dia}.)

For each edge of $P$ joining one cell of $\Phi(b')$ and one cell of $\Phi(b'')$, with $b' \neq b''$, take the edge $b'b''$ of $H$. Since $P$ visits $\Phi(b)$ for all $b$, these edges form a Hamiltonian path in~$H$. Denote this path by $\varrho$.

From this point on, our plan for the proof will be as follows: First we obtain a complete description of the structure of $\varrho$. We do this by means of a series of mostly local considerations, starting on the boundary of $H$ and then working our way in. Once we are done, we determine what paths $P$ in the original graph $G$ are associated with each path $\varrho$.

We continue with the details.

We define the $i$-th \emph{frame} of $B$, denoted $F_i$, to be the subset \[F_i = [i; k - i - 1]^2 \setminus [i + 1; k - i - 2]^2\] of $B$. Thus $F_0$, $F_1$, $\ldots$, $F_{\lceil k/2 \rceil - 1}$ form a partitioning of $B$.

We denote the four corner cells of $F_i$ by $a_i = (i, i)$, $b_i = (k - i - 1, i)$, $c_i = (k - i - 1, k - i - 1)$, and $d_i = (i, k - i - 1)$.

\begin{figure}[ht] \null \hfill \begin{subfigure}{40pt} \centering \includegraphics{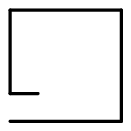} \caption{} \label{snail-a} \end{subfigure} \hfill \begin{subfigure}{40pt} \centering \includegraphics{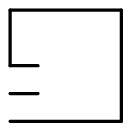} \caption{} \label{snail-b} \end{subfigure} \hfill \begin{subfigure}{40pt} \centering \includegraphics{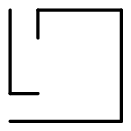} \caption{} \label{snail-c} \end{subfigure} \hfill \begin{subfigure}{40pt} \centering \includegraphics{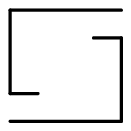} \caption{} \label{snail-d} \end{subfigure} \hfill \begin{subfigure}{40pt} \centering \includegraphics{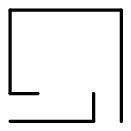} \caption{} \label{snail-e} \end{subfigure} \hfill \null \caption{} \label{snails} \end{figure}

We say that $F_i$ is of type I, II, III, IV, or V (relative to $\varrho$) when $\varrho$ contains the following subpaths:

For type I, $a_i \spring b_i \spring c_i \spring d_i \spring (a_i + (0, 1))$---$(a_i + (1, 1))$. (Figure \ref{snail-a}.)

For type II, $a_i \spring b_i \spring c_i \spring d_i \spring (a_i + (0, 2))$---$(a_i + (1, 2))$ and $(a_i + (0, 1))$---$(a_i + (1, 1))$. (Figure \ref{snail-b}.)

For type III, $a_i \spring b_i \spring c_i \spring (d_i + (1, 0))$---$(d_i + (1, -1))$ and $d_i \spring (a_i + (0, 1))$---$(a_i + (1, 1))$. (Figure \ref{snail-c}.)

For type IV, $a_i \spring b_i \spring (c_i + (0, -1))$---$(c_i + (-1, -1))$ and $c_i \spring d_i \spring (a_i + (0, 1))$---$(a_i + (1, 1))$. (Figure \ref{snail-d}.)

For type V, $a_i \spring (b_i + (-1, 0))$---$(b_i + (-1, 1))$ and $b_i \spring c_i \spring d_i \spring (a_i + (0, 1))$---$(a_i + (1, 1))$. (Figure \ref{snail-e}.)

We say that $\varrho$ itself is of type I when all of $F_0$, $F_1$, $\ldots$, $F_{\lceil k/2 \rceil - 2}$ are of type I. Furthermore, let $\mathfrak{T}$ be one of the symbols II, III, IV, and V, and let $s$ be a nonnegative integer with $0 \le s \le \lceil k/2 \rceil - 2$. We say that $\varrho$ is of type $\mathfrak{T}(s)$ when all of $F_0$, $F_1$, $\ldots$, $F_{s - 1}$ are of type I and all of $F_s$, $F_{s + 1}$, $\ldots$, $F_{\lceil k/2 \rceil - 2}$ are of type $\mathfrak{T}$. When $\varrho$ is of one of the $4\lceil k/2 \rceil - 3$ types we have just listed, we say that it is \emph{regular}.

\begin{figure}[ht] \centering \includegraphics{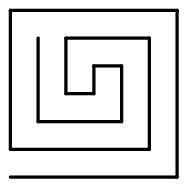} \caption{} \label{spiral-7x7} \end{figure}

For example, Figure \ref{spiral-7x7} shows the unique $\varrho$ of type $\operatorname{III}(1)$ when $k = 7$.

Observe that, when $k$ is even, $\varrho$ cannot be of type $\operatorname{IV}(s)$, for any $s$, as $F_{\lceil k/2 \rceil - 2}$ being of type IV prevents $\varrho$ from visiting all cells of $F_{\lceil k/2 \rceil - 1}$. Similarly, when $k$ is odd, $\varrho$ cannot be of type $\operatorname{II}(s)$, for any $s$, as the frame $F_{\lceil k/2 \rceil - 2}$ is too small to be of type II.

In all other cases, if $k$ and the type of $\varrho$ are fixed, there is a unique $\varrho$ of that type, with one exception: When $k$ is even, there are two paths $\varrho$ of type I, differing by just one edge within the innermost frame $F_{\lceil k/2 \rceil - 1}$.

\medskip

\begin{lemma} \label{turn} Suppose that $\varrho$ makes a turn at cell $b$ of $B$, for concreteness by means of $(b + (0, 1))$---$b$---$(b + (1, 0))$. Then the two cells of $\Phi(b)$ in $P$ are $2b + (0, 1)$ and $2b + (1, 0)$. \end{lemma}

\medskip

Of course, similar claims hold for the other three possible turns at $b$ as well.

\medskip

\begin{proof} Suppose, for the sake of contradiction, that $2b + (1, 1) \in P$. Since $(b + (0, 1))$---$b$---$(b + (1, 0)) \subseteq \varrho$, we get that there are two cells $a'$ and $a''$ with $a' \in \Phi(b + (0, 1))$, $a''\in \Phi(b + (1, 0))$, and $a'$---$(2b + (1, 1))$---$a'' \subseteq P$. However, then $P$ cannot visit any cells of $\Phi(b)$ other than $2b + (1, 1)$, a contradiction.

\begin{figure}[ht] \centering \includegraphics{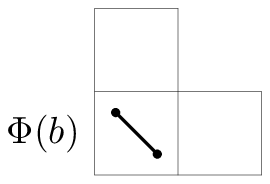} \caption{} \label{smooth-turn} \end{figure}

Thus $2b + (1, 1) \not \in P$. Since $(b + (0, 1))$---$b$---$(b + (1, 0)) \subseteq \varrho$, it follows that $2b + (0, 1) \in P$ and $2b + (1, 0) \in P$, as needed. (Figure \ref{smooth-turn}.) \end{proof}

\medskip

\begin{lemma} \label{spiral} A symmetry of $B$ maps $\varrho$ onto a regular Hamiltonian path in $H$. \end{lemma}

\medskip

\begin{proof} First we consider the outermost frame $F_0$ of $B$. Observe that all edges of $H$ within $F_0$ form a cycle $E$.

Suppose, for the sake of contradiction, that $(0, y)$---$(0, y + 1) \not \in \varrho$ for some $y$ with $2 \le y \le k - 4$.

\begin{figure}[ht] \centering \includegraphics{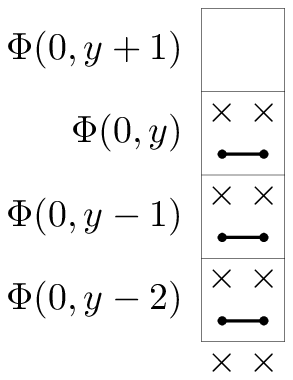} \caption{} \label{three-endpoints} \end{figure}

Then $P$ must miss both cells in at least one of the two pairs $\{(0, 2y + 1), (1, 2y + 1)\}$ and $\{(0, 2y + 2), (1, 2y + 2)\}$. Suppose, for concreteness, that this is true of the former pair; the other case is similar. (Figure \ref{three-endpoints}.)

Since $P$ visits two cells in $\Phi(0, y)$, we get that $(0, 2y) \in P$ and $(1, 2y) \in P$. So $(0, 2y - 1) \not \in P$ and $(1, 2y - 1) \not \in P$. Iterating this sequence of observations, we see that all of $(0, 2y - 2)$, $(1, 2y - 2)$, $(0, 2y - 4)$, and $(1, 2y - 4)$ are in $P$ while all of $(0, 2y - 3)$, $(1, 2y - 3)$, $(0, 2y - 5)$, and $(1, 2y - 5)$ are not. (When $y = 2$, the lattermost couple of cells will lie outside of $A$. Of course, then they will be outside of $P$, too.)

Thus all three of $(0, 2y)$, $(0, 2y - 2)$, and $(0, 2y - 4)$ must be endpoints of $P$, a contradiction.

Consequently, $(0, y)$---$(0, y + 1) \in \varrho$ for all $y$ with $2 \le y \le k - 4$. By symmetry, it follows that $\varrho$ contains all edges of $E$ except for, possibly, $(0, 0)$---$(0, 1)$, $(0, 1)$---$(0, 2)$, and their images under the symmetries of $B$.

On the other hand, at least one edge of $E$ must be outside of $\varrho$.

If $(0, 0)$---$(0, 1) \not \in \varrho$, then $(0, 0)$ is an endpoint of $\varrho$ and $\Phi(0, 0)$ contains an endpoint of $P$.

If $(0, 1)$---$(0, 2) \not \in \varrho$, then $P$ must miss both cells in at least one of the two pairs $\{(0, 3), (1, 3)\}$ and $\{(0, 4), (1, 4)\}$. If the latter, then with $k \ge 5$ we arrive at a contradiction as before. If the former, then the same reasoning as before shows that $\Phi(0, 0)$ and $\Phi(0, 1)$ must each contain an endpoint of $P$. When $k = 4$, the same conclusion holds up to reflection with respect to the horizontal axis of symmetry of $A$.

Of course, these observations apply also to all images of the edges $(0, 0)$---$(0, 1)$ and $(0, 1)$---$(0, 2)$ under the symmetries of $B$.

Since there are only two endpoints of $P$, it follows that $\varrho$ cannot omit too many edges of $E$. We are left to consider the following cases, up to the symmetries of $B$:

\smallskip

\emph{Case 1}. The only edge of $E$ outside of $\varrho$ is $(0, 0)$---$(0, 1)$. Then $F_0$ is of type I.

\smallskip

\emph{Case 2}. The only edge of $E$ outside of $\varrho$ is $(0, 1)$---$(0, 2)$, and $k \ge 4$. Then $F_0$ is of type II.

\smallskip

\emph{Case 3}. The only edges of $E$ outside of $\varrho$ are $(0, 0)$---$(0, 1)$ and $(0, k - 1)$---$(1, k - 1)$. Then $F_0$ is of type III.

\smallskip

\emph{Case 4}. The only edges of $E$ outside of $\varrho$ are $(0, 0)$---$(0, 1)$ and $(k - 1, k - 2)$---$(k - 1, k - 1)$. Then $F_0$ is of type IV.

\smallskip

\emph{Case 5}. The only edges of $E$ outside of $\varrho$ are $(0, 0)$---$(0, 1)$ and $(k - 2, 0)$---$(k - 1, 0)$. Then $F_0$ is of type V.

\smallskip

\emph{Cases 6--9}. The only edges of $E$ outside of $\varrho$ are $(0, 0)$---$(0, 1)$ and one of $(0, 0)$---$(1, 0)$, $(0, k - 2)$---$(0, k - 1)$, $(k - 2, k - 1)$---$(k - 1, k - 1)$, and $(k - 1, 0)$---$(k - 1, 1)$. The first one of these cases cannot occur because then cell $(0, 0)$ becomes isolated in $\varrho$. The other three cannot occur because, in each one of them, part of the edges of $E$ in $\varrho$ form a subpath of $\varrho$ which contains both endpoints of $\varrho$ but does not coincide with $\varrho$.

\smallskip

With this, we have established that $F_0$ is of one of our types relative to $\varrho$, up to the symmetries of $B$.

We finish the proof by induction on $k$. Our base cases are $k = 3$ and $k = 4$, when there is nothing left to prove. For the induction step, suppose that $k \ge 5$ and that we have already settled the question on all smaller boards.

Let $A^\star$ be the concentric subboard of $A$ of size $(n - 4) \times (n - 4)$ given by $A^\star = [2; n - 3]^2$ and let $G^\star$ be the king graph on $A^\star$.

When $F_0$ is of type I, the cells of $P$ outside of $A^\star$ form a subpath of $P$ containing one endpoint of $P$. When $F_0$ is of one of the remaining four types, the cells of $P$ outside of $A^\star$ form two subpaths of $P$ containing the two endpoints of $P$. Either way, we see that the restriction $P^\star$ of $P$ to $A^\star$ is a subpath of $P$. Since $P$ is a snake path in $G$, also $P^\star$ is a snake path in $G^\star$. Furthermore, $P^\star$ will be of the greatest possible length within $G^\star$.

Let $B^\star$ be the concentric subboard of $B$ of size $(k - 2) \times (k - 2)$ given by $B^\star = [1; k - 2]^2 = B \setminus F_0$ and let $H^\star$ be the grid graph on $B^\star$. By the same reasoning as above, the restriction $\varrho^\star$ of $\varrho$ to $B^\star$ is a subpath of $\varrho$. Furthermore, $\varrho^\star$ is a Hamiltonian path in $H^\star$, and $\varrho^\star$ and $P^\star$ are related in the same way as $\varrho$ and $P$, in the sense that $\varrho^\star$ visits the cells of $B^\star$ in the same order as $P^\star$ visits their corresponding $2 \times 2$ subboards of $A^\star$.

By the induction hypothesis, either $\varrho^\star$ or an image of it under a symmetry of $B^\star$ must be regular relative to $B^\star$.

We consider five cases for the type of $F_0$.

\smallskip

\emph{Case 1}. $F_0$ is of type I.

Then $(0, 1)$---$(1, 1) \in \varrho$ and $(1, 1)$ is an endpoint of $\varrho^\star$.

Suppose, for the sake of contradiction, that $(1, 1)$---$(1, 2) \in \varrho$ and $\varrho$ makes a turn at $(1, 1)$.

If $\varrho$ also makes a turn at $(2, 1)$, then by Lemma \ref{turn} we get that $(3, 3) \in P$ and $(4, 3) \in P$, in contradiction with $(1, 1)$---$(2, 1) \not \in \varrho$.

Otherwise, if $\varrho$ does not make a turn at $(2, 1)$, then $(2, 1)$ is an endpoint of $\varrho^\star$. By the induction hypothesis, it follows that $F_1$ satisfies the reflection, with respect to the line $x = y$, of the conditions defining type II. Thus $\varrho$ makes a turn at $(3, 1)$, and by Lemma \ref{turn} we get that $(2, 2)$---$(3, 3)$---$a' \subseteq P$ and $(7, 2)$---$(6, 3)$---$a'' \subseteq P$ for two cells $a'$ and $a''$ with $a' \in \Phi(1, 2)$ and $a'' \in \Phi(3, 2)$. Consequently, $P$ cannot visit any cells in $\Phi(2, 1)$, another contradiction.

We conclude that $(1, 1)$---$(1, 2) \not \in \varrho$ and $(1, 1)$---$(2, 1) \in \varrho$. By the induction hypothesis, it follows that $\varrho^\star$ is regular. (As in, it is $\varrho^\star$ itself that is regular, rather than some image of it under a symmetry of $B^\star$.) Then $\varrho$ is regular as well, and of the same type as~$\varrho^\star$.

\smallskip

\emph{Case 2}. $F_0$ is of type II.

Then $(0, 1)$---$(1, 1) \in \varrho$, $(0, 2)$---$(1, 2) \in \varrho$, and $(1, 1)$ and $(1, 2)$ are the two endpoints of $\varrho^\star$. By the induction hypothesis, it follows that $\varrho^\star$ is of type $\operatorname{II}(0)$, and so is $\varrho$ as well.

\smallskip

\emph{Case 3}. $F_0$ is of type III.

Then $(0, 1)$---$(1, 1) \in \varrho$, $(1, k - 1)$---$(1, k - 2) \in \varrho$, and $(1, 1)$ and $(1, k - 2)$ are the two endpoints of $\varrho^\star$. By the induction hypothesis, it follows that either $\varrho^\star$ is of type $\operatorname{III}(0)$, or its reflection with respect to the horizontal axis of symmetry of $B^\star$ is. (Here, we take into account the fact that types $\operatorname{III}(0)$ and $\operatorname{V}(0)$ are related by quarter-turn rotation.) If the latter, then $\varrho$ makes a turn at both cells $(1, 1)$ and $(2, 1)$, and we get a contradiction as in Case 1. Thus $\varrho^\star$ is of type $\operatorname{III}(0)$, and so is $\varrho$ as well.

\smallskip

\emph{Case 4}. $F_0$ is of type IV.

Then $(0, 1)$---$(1, 1) \in \varrho$, $(k - 1, k - 2)$---$(k - 2, k - 2) \in \varrho$, and $(1, 1)$ and $(k - 2, k - 2)$ are the two endpoints of $\varrho^\star$. By the induction hypothesis, it follows that either $\varrho^\star$ is of type $\operatorname{IV}(0)$, or its reflection with respect to the line $x = y$ is. From here, the analysis continues as in Case 3, and in the end we obtain that both of $\varrho^\star$ and $\varrho$ are of type $\operatorname{IV}(0)$.

\smallskip

\emph{Case 5}. $F_0$ is of type V.

This case is analogous to Case 3, and in it both of $\varrho^\star$ and $\varrho$ are of type $\operatorname{V}(0)$.

\smallskip

With this, the induction step is complete. \end{proof}

\medskip

\begin{lemma} \label{blowup} Let $\varrho$ be a regular Hamiltonian path in $H$. Then the number of longest snake paths $P$ in $G$ associated with $\varrho$ is as shown in Table \ref{count}.

\begin{table}[ht] \centering \begin{tabular}{c|c|c|} type of $\varrho$ & $k$ even & $k$ odd\\ \hline I & $8 + 1$ & $8$\\ \hline $\operatorname{II}(s)$ & $6$ & ---\\ \hline $\operatorname{III}(s)$ & $1$ & $1$\\ \hline $\operatorname{IV}(s)$ & --- & $6$\\ \hline $\operatorname{V}(s)$ & $1$ & $1$\\ \hline \end{tabular} \caption{} \label{count} \end{table}

These paths are pairwise distinct under the symmetries of $A$, with the following exceptions: The unique $P$ with $\varrho$ of type $\operatorname{III}(0)$ and the unique $P$ with $\varrho$ of type $\operatorname{V}(0)$ are related by quarter-turn rotation; and, when $k$ is odd, two pairs of paths $P$ with $\varrho$ of type $\operatorname{IV}(0)$ are related by central symmetry. \end{lemma}

\medskip

(The entry of Table \ref{count} for $k$ even and $\varrho$ of type I includes one summand for each regular path $\varrho$ of that type. For the parametrised types, the total count does not depend on the value of the parameter.)

\medskip

\begin{proof} First we strengthen the claim as follows: Let $b$ be an arbitrary cell of $F_0$. Then, for all $k$ with $k \ge 5$:

(i) When $\varrho$ goes straight through $b$, the two cells of $\Phi(b)$ in $P$ are the ones on the boundary of $A$, in the union of row $0$, row $n - 1$, column $0$, and column $n - 1$ of $A$; and

(ii) When $b$ is an endpoint of $\varrho$, by the definitions of our types we get that $b$ is one of $(0, 0)$, $(0, 1)$, $(k - 1, 0)$, $(k - 1, k - 1)$, and $(0, k - 1)$. Then the two cells of $\Phi(b)$ in $P$ are, respectively, the ones in row $0$, row $2$, column $n - 1$, row $n - 1$, and column $0$ of $A$.

For the proof, we work by induction on $k$.

\begin{figure}[ht] \null \hfill \begin{subfigure}{85pt} \centering \includegraphics{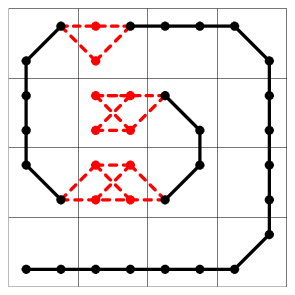} \caption{} \label{base-case-even-a-i} \end{subfigure} \hfill \begin{subfigure}{85pt} \centering \includegraphics{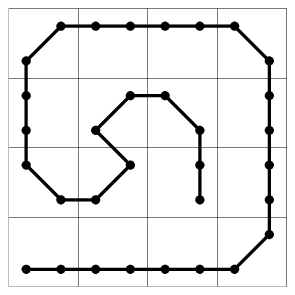} \caption{} \label{base-case-even-a-ii} \end{subfigure} \hfill \null\\ \vspace{\baselineskip}\\ \null \hfill \begin{subfigure}{85pt} \centering \includegraphics{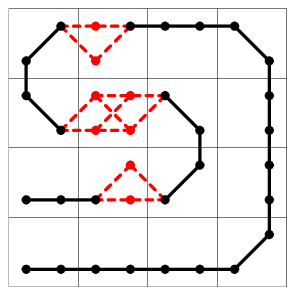} \caption{} \label{base-case-even-b} \end{subfigure} \hfill \begin{subfigure}{85pt} \centering \includegraphics{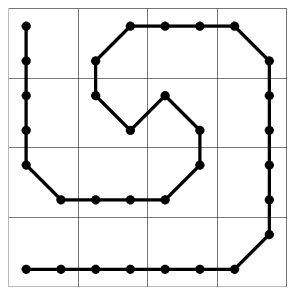} \caption{} \label{base-case-even-c} \end{subfigure} \hfill \begin{subfigure}{85pt} \centering \includegraphics{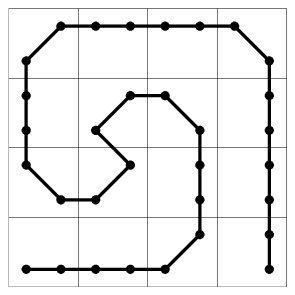} \caption{} \label{base-case-even-e} \end{subfigure} \hfill \null \caption{} \label{base-case-even} \end{figure}

Our base cases are $k = 3$ and $k = 4$, when the strengthening is irrelevant and the original claim follows by direct examination of all cases. These are shown in Figures \ref{base-case-even} and \ref{base-case-odd}. For each situation, we depict the cells and edges which are forced to belong to $P$ in black, and the optional cells and edges of $P$ in red. Figures \ref{base-case-even-a-i}--\subref{base-case-even-e} correspond to the two paths $\varrho$ of type I and the paths $\varrho$ of types $\operatorname{II}(0)$, $\operatorname{III}(0)$, and $\operatorname{V}(0)$, respectively, while Figures \ref{base-case-odd-a}--\subref{base-case-odd-e} correspond to the paths $\varrho$ of types I, $\operatorname{III}(0)$, $\operatorname{IV}(0)$, and $\operatorname{V}(0)$, respectively.

\begin{figure}[ht] \null \hfill \begin{subfigure}{65pt} \centering \includegraphics{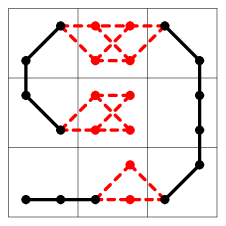} \caption{} \label{base-case-odd-a} \end{subfigure} \hfill \begin{subfigure}{65pt} \centering \includegraphics{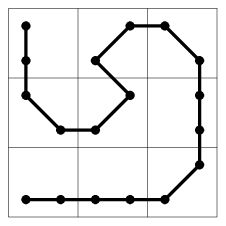} \caption{} \label{base-case-odd-c} \end{subfigure} \hfill \begin{subfigure}{65pt} \centering \includegraphics{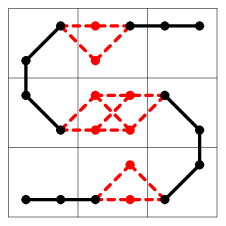} \caption{} \label{base-case-odd-d} \end{subfigure} \hfill \begin{subfigure}{65pt} \centering \includegraphics{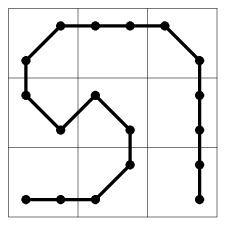} \caption{} \label{base-case-odd-e} \end{subfigure} \hfill \null \caption{} \label{base-case-odd} \end{figure}

For the induction step, suppose that $k \ge 5$ and that we have already settled the question on all smaller boards. Define $A^\star$, $B^\star$, $G^\star$, $H^\star$, $P^\star$, and $\varrho^\star$ as in the proof of Lemma \ref{spiral}.

We begin with parts (i) and (ii) of the strengthening.

For (i), consider a cell $b$ of $F_0$ such that $\varrho$ goes straight through $b$, using two edges on $H$ within $F_0$. Let $c$ be the unique neighbour of $b$ in $H$ which belongs to $F_1$, let $b'$ and $b''$ be the two cells of $\Phi(b)$ which are adjacent in $G$ to a cell in $\Phi(c)$, and, conversely, let $c'$ and $c''$ be the two cells of $\Phi(c)$ which are adjacent in $G$ to a cell in $\Phi(b)$.

Then at least one of $c'$ and $c''$ must be in $P$. Indeed, when $\varrho^\star$ makes a turn at $c$, this follows by Lemma \ref{turn}. Otherwise, when $\varrho^\star$ goes straight through $c$ using two edges of $H$ within $F_1$, or when $c$ is an endpoint of $\varrho^\star$, it follows by the catalogue in Figures \ref{base-case-even} and \ref{base-case-odd} when $k = 5$ or $k = 6$, and by the induction hypothesis when $k \ge 7$.

Since $bc \not \in \varrho$ and at least one of $c'$ and $c''$ is in $P$, we conclude that $b' \not \in P$ and $b'' \not \in P$. So the other two cells of $\Phi(b)$ must be in $P$, confirming (i).

For (ii), consider a cell $b$ of $F_0$ such that $b$ is an endpoint of $\varrho$. When $b$ is $(0, 0)$, $(0, 1)$, $(k - 1, 0)$, $(k - 1, k - 1)$, or $(0, k - 1)$, define $c$ to be $(0, 1)$, $(0, 2)$, $(k - 2, 0)$, $(k - 1, k - 2)$, or $(1, k - 1)$, respectively, and also define cells $b'$, $b''$, $c'$, and $c''$ relative to $b$ and $c$ as in our treatment of (i).

By the definitions of our types, $bc \not \in \varrho$. Furthermore, with one exception, to be considered shortly, $\varrho$ makes a turn at $c$. By Lemma \ref{turn}, it follows that exactly one of $c'$ and $c''$ is in $P$. From here, as in our treatment of (i), $bc \not \in \varrho$ implies $b' \not \in P$ and $b'' \not \in P$, and so the other two cells of $\Phi(b)$ must be in $P$, as needed. The unique exception is when $F_0$ is of type II, $b = (0, 0)$, and $c = (0, 1)$. However, then $c' \in P$ and $c'' \in P$ by what we just proved applied to the endpoint $(0, 1)$ of $\varrho$, and once again (ii) is confirmed.

With this, we have established both parts (i) and (ii) of the strengthening.

For the original claim, observe that the type of $\varrho$ determines the type of $F_0$ uniquely. We are left to show that the type of $F_0$, and the subpath $P^\star$, determine $P$ uniquely. We only need to look at the cells of $P$ in $A \setminus A^\star$, that is, in the union of the subboards $\Phi(b)$ of $A$ with $b \in F_0$.

Let, then, $b$ be an arbitrary cell of $F_0$. When $\varrho$ makes a turn at $b$, the cells of $P$ in $\Phi(b)$ are uniquely determined by Lemma \ref{turn}. Otherwise, when $\varrho$ goes straight through $b$, or when $b$ is an endpoint of $\varrho$, the desired uniqueness follows by parts (i) and (ii) of the strengthening, respectively. The induction step is complete. \end{proof}

\medskip

\begin{proof*}{Proof of Theorem \ref{king-even-path}} We get the number of essentially distinct paths by Lemmas \ref{spiral} and \ref{blowup}, summing over all possible types of $\varrho$ and then subtracting out the duplicates specified in the statement of Lemma \ref{blowup}.

To convert this into the total number of paths, we analyse symmetries.

When $k$ is even, all of the $P$'s are asymmetric.

Otherwise, when $k$ is odd, two of the $P$'s which correspond to the unique $\varrho$ of type $\operatorname{IV}(0)$ are preserved under central symmetry, but not under any other symmetries of $A$; and all other $P$'s corresponding to regular $\varrho$'s are asymmetric. \end{proof*}

\medskip

Note that this argument also confirms our remark in the introduction regarding the number of essentially distinct longest snake paths in $G$.

The proof of Theorem \ref{king-even-path} shows that, roughly speaking, every longest snake path $P$ in $G$ behaves as follows: Outside of some concentric square subboard of $A$, it is shaped as a simple spiral; and then, within that subboard, it is shaped as a double spiral instead. One endpoint of $P$ lies on the boundary of $A$, and the other one lies on the boundary between the two spirals.

%% file: snake-04-king-odd-i.tex
\section{King Graphs on Odd Boards I} \label{king-odd-i}

In this section, we prove Theorem \ref{king-odd-path}. Let $n$ be an odd positive integer with $n = 2k - 1$, let $A$ be the standard square board of side $n$, and let $G$ be the king graph on $A$.

\medskip

\begin{proof*}{Proof of the lower bound for Theorem \ref{king-odd-path}} Consider the set of all cells of $A$ of the form $(x, y)$ with $1 \le x \le n - 2$ and $y$ even. The king graph on it is the disjoint union of $k$ paths. To obtain the vertex set of a single snake path in $G$, we add in also the cells $(0, 0)$; $(0, n - 1)$ if $k$ is even and $(n - 1, n - 1)$ if $k$ is odd; $(0, y)$ in $A$ with $y \equiv 3 \pmod 4$; and $(n - 1, y)$ in $A$ with $y \equiv 1 \pmod 4$. \end{proof*}

\medskip

\begin{figure}[ht] \centering \includegraphics{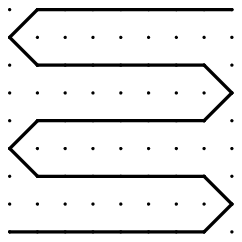} \caption{} \label{king-9x9} \end{figure}

For example, Figure \ref{king-9x9} shows $n = 9$.

We obtained the upper bound of Theorem \ref{king-even-path} by summing over some subsets of $A$ such that, for every snake path $P$, the part of $P$ within each subset must be small. Our approach to the upper bound of Theorem \ref{king-odd-path} will follow a similar strategy, albeit with significant complications. Instead of subsets of $A$, we sum over subgraphs of $G$. Furthermore, our notion of smallness will be somewhat unusual: We consider the total number of certain cells and edges of $P$ within each subgraph.

We call a cell $(x, y)$ of $A$ \emph{even} when both of $x$ and $y$ are even, and \emph{odd} when both of $x$ and $y$ are odd. We also call an edge of $G$ \emph{regular} when it is not incident with an odd cell.

\begin{figure}[ht] \null \hfill \begin{subfigure}{20pt} \centering \includegraphics{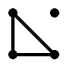} \caption{} \label{block-little} \end{subfigure} \hfill \begin{subfigure}{30pt} \centering \includegraphics{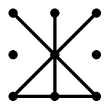} \caption{} \label{block-large} \end{subfigure} \hfill \null \caption{} \label{blocks-figure} \end{figure}

Given an even cell $a = (z, z)$ of $A$ with $z \le k - 2$, we write $\nosovka(a)$ for the subgraph of $G$ with vertices $a + [0; 1]^2$ whose edges join $a$, $a + (0, 1)$, and $a + (1, 0)$ pairwise. Thus $\nosovka(a)$ contains four cells, one of which is odd, and three edges, all of which are regular. (Figure \ref{block-little}.)

For every symmetry $\pi$ of $A$, if $b = \pi(a)$, then we also define $\nosovka(b) = \pi(\nosovka(a))$. We call each subgraph of $G$ of this form a \emph{little block}.

Given an even cell $a = (x, y)$ of $A$ with $x > y$ and $x + y \le n - 3$, we write $\nosovka(a)$ for the subgraph of $G$ with vertices $a + [-1; 1] \times [0; 2]$ whose edges join $a$ to $a + (-1, 0)$ and $a + (1, 0)$ as well as $a + (0, 1)$ to all elements of the set $a + [-1; 1] \times \{0, 2\}$. Thus $\nosovka(a)$ contains nine cells, two of which are odd, and eight edges, all of which are regular. (Figure \ref{block-large}.)

For every symmetry $\pi$ of $A$, if $b = \pi(a)$, then we also define $\nosovka(b) = \pi(\nosovka(a))$. We call each subgraph of $G$ of this form a \emph{large block}.

\begin{figure}[ht] \null \hfill \begin{subfigure}{125pt} \centering \includegraphics{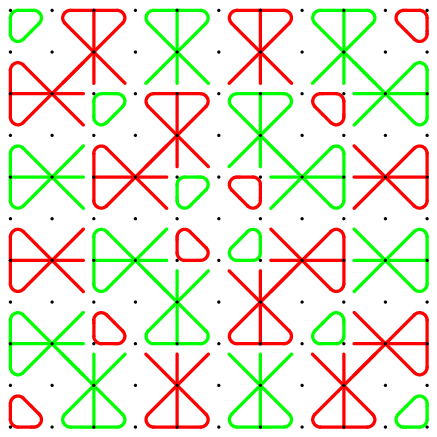} \caption{} \label{blocks-11x11} \end{subfigure} \hfill \begin{subfigure}{150pt} \centering \includegraphics{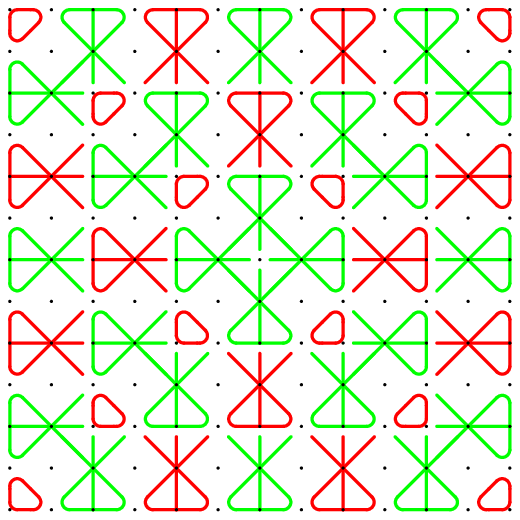} \caption{} \label{blocks-13x13} \end{subfigure} \hfill \null \caption{} \label{blocks-examples} \end{figure}

For example, Figure \ref{blocks-examples} shows $n = 11$ and $n = 13$. (The colouring is only for clarity.)

Observe that every odd cell belongs to at least two blocks and every regular edge belongs to at least one block.

Consider an arbitrary snake path $P$ in $G$. For every block $\nosovka$, let $w_\text{Cell}(\nosovka)$ be the number of odd cells of $P$ in $\nosovka$, let $w_\text{Edge}(\nosovka)$ be the number of regular edges of $P$ in $\nosovka$, and let $w(\nosovka) = w_\text{Cell}(\nosovka) + w_\text{Edge}(\nosovka)$.

\medskip

\begin{lemma} \label{blocks} When $\nosovka$ is a little block, $w(\nosovka) \le 1$. Otherwise, when $\nosovka$ is a large block, $w(\nosovka) \le 2$. \end{lemma}

\medskip

\begin{proof} By direct examination of all cases. \end{proof}

\medskip

\begin{proof*}{Proof of the upper bound for Theorem \ref{king-odd-path}} Let $w_\text{Cell}(P)$ be the number of odd cells in $P$ and let $w_\text{Edge}(P)$ be the number of regular edges in $P$. Then the length of $P$ is bounded from above by $2w_\text{Cell}(P) + w_\text{Edge}(P)$. We sum the inequalities of Lemma \ref{blocks} over all blocks $\nosovka$, and we obtain that the latter expression cannot exceed $(n^2 - 1)/2$. \end{proof*}

%% file: snake-05-king-odd-ii.tex
\section{King Graphs on Odd Boards II} \label{king-odd-ii}

In this section, we give a complete description of the longest snake paths in king graphs on odd square boards. We use the same notations as in Section \ref{king-odd-i}.

We begin with brief overviews of two relevant topics.

Let $s$ be a positive integer and let $\sigma$ be a permutation of $[0; s - 1]$. For each $i$ with $0 \le i \le s - 2$, draw a semicircle with endpoints $(0, \sigma^{-1}(i))$ and $(0, \sigma^{-1}(i + 1))$ which lies on the left of the coordinate axis $Oy$ when $i$ is even and on its right otherwise, when $i$ is odd. The union of all such semicircles is a curve $\kappa$ in the plane with endpoints $(0, \sigma^{-1}(0))$ and $(0, \sigma^{-1}(s - 1))$. When this curve does not intersect itself, $\sigma$ is a \emph{stamp-folding permutation}.

\begin{figure}[ht] \null \hfill \begin{minipage}{50pt} \centering \includegraphics{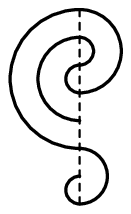} \caption{} \label{stamp} \end{minipage} \hfill \begin{minipage}{90pt} \centering \includegraphics{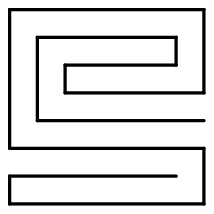} \caption{} \label{fewest-turns} \end{minipage} \hfill \null \end{figure}

For example, Figure \ref{stamp} shows $\kappa$ when $\sigma$ is the permutation $1$, $0$, $2$, $7$, $4$, $5$, $6$, $3$.

Combinatorially, $\sigma$ is a stamp-folding permutation if and only if there are no $i$ and $j$ of the same parity with $0 \le i \le s - 2$ and $0 \le j \le s - 2$ such that exactly one of $i$ and $i + 1$ lies between $j$ and $j + 1$ in $\sigma$.

The intuition is as follows: Imagine a paper strip of size $1 \times s$ formed out of $s$ stamps of size $1 \times 1$ each. We fold this strip along the perforations between stamps so that all stamps come to lie on top of one another. Then $\sigma$ is a stamp-folding permutation if and only if it can be obtained as the permutation of the stamps within the resulting stack. The points of $\kappa$ on the coordinate axis $Oy$ correspond to the $s$ stamps, and the semicircles of $\kappa$ correspond to the $s - 1$ creases between stamps.

Stamp-folding permutations have been studied extensively.

Consider a Hamiltonian path in the grid graph $\Gamma$ of size $s \times s$. The smallest number of turns that such a path can make is $2s - 2$. \cite{J} We proceed to review some properties of the paths which attain this minimum. \cite{G}

(The author noticed the connection between stamp-folding permutations and fewest-turn Hamiltonian paths when he obtained Theorems \ref{king-odd-path} and \ref{king-odd-structure}. Independently, \cite{G} was published before the present work was written.)

Let $\alpha$ be a fewest-turn Hamiltonian path in $\Gamma$. We partition $\alpha$ into $2s - 1$ subpaths at the $2s - 2$ cells where it makes a turn. (Each cell with a turn belongs to two subpaths.) We call these subpaths the \emph{segments} of $\alpha$. Thus in each segment of $\alpha$ either all edges are horizontal or all edges are vertical, and segments of these two types alternate.

Let us call $\alpha$ \emph{mostly-horizontal} when it consists of $s$ horizontal segments and ${s - 1}$ vertical segments, and \emph{mostly-vertical} otherwise, when it is the other way around. Suppose, for concreteness, that $\alpha$ is mostly-horizontal.

Then every row of the board contains exactly one horizontal segment of $\alpha$. Orient $\alpha$ arbitrarily, and then number its segments from $0$ to $s - 1$ in the order in which they occur along $\alpha$. When we assign to each row of the board the number of its horizontal segment of $\alpha$, we obtain a stamp-folding permutation. Conversely, every stamp-folding permutation corresponds in this way to exactly two oriented mostly-horizontal fewest-turn Hamiltonian paths. The two are reflections of one another with respect to the vertical axis of symmetry of the board.

Explicitly, the correspondence is as follows: For each element $i$ of $[0; s - 1]$, let $\omega_\text{Left}(i)$ be the number of even nonnegative integers $j$ such that $i$ lies between $j$ and $j + 1$ in $\sigma$, and define $\omega_\text{Right}(i)$ similarly, but with $j$ odd. In one of the two oriented paths associated with $\sigma$, the path's $i$-th horizontal segment goes from $(s - \omega_\text{Right}(i) - 1, \sigma^{-1}(i))$ to $(\omega_\text{Left}(i), \sigma^{-1}(i))$ for all even $i$, and it goes in the opposite direction for all odd $i$. The second oriented path associated with $\sigma$ can be obtained by reflection.

For example, Figure \ref{fewest-turns} shows this for the stamp-folding permutation of Figure \ref{stamp}.

One corollary of this connection is that, for all $s$ with $s \ge 2$, the number of fewest-turn Hamiltonian paths in the grid graph of size $s \times s$ is twice the number of stamp-folding permutations of $s$ elements.

This concludes the two overviews, and we return to our main topic.

Let $H$ be the grid graph of size $k \times k$ and let $\varrho$ be a Hamiltonian path in it. (Recall from Section \ref{king-odd-i} that $k = \lceil n/2 \rceil$.)

We say that the cycle $a'a''b'b''$ in $H$ is \emph{free} (relative to $\varrho$) when $a'a''$ is an edge of $\varrho$; the other three edges of the cycle are outside of $\varrho$; both of $a'$ and $a''$ are cells that $\varrho$ goes straight through (so that, if $a' = (c' + a'')/2$ and $a'' = (a' + c'')/2$, then $c' \spring c'' \subseteq \varrho$); and both of $b'$ and $b''$ are cells where $\varrho$ makes a turn.

\begin{figure}[t] \null \hfill \begin{subfigure}{85pt} \centering \includegraphics{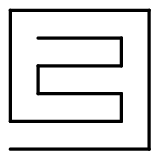} \caption{} \label{rho} \end{subfigure} \hfill \begin{subfigure}{85pt} \centering \includegraphics{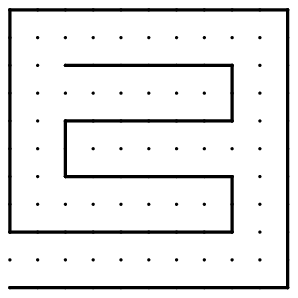} \caption{} \label{phi} \end{subfigure} \hfill \null\\ \vspace{\baselineskip}\\ \null \hfill \begin{subfigure}{85pt} \centering \includegraphics{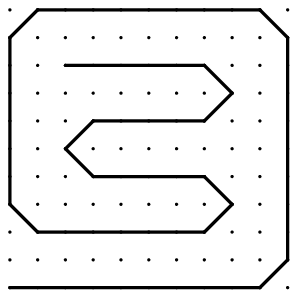} \caption{} \label{psi} \end{subfigure} \hfill \begin{subfigure}{85pt} \centering \includegraphics{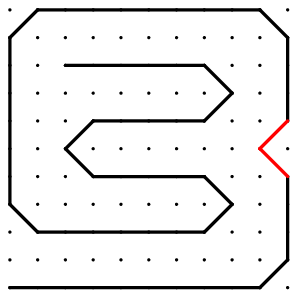} \caption{} \label{lift} \end{subfigure} \hfill \null \caption{} \label{lift-examples} \end{figure}

We proceed to associate $\varrho$ with certain paths in $G$. Figure \ref{rho} shows one example of a Hamiltonian path $\varrho$ in $H$, and Figures \ref{phi}--\subref{lift} track the series of definitions given below.

For each edge $a'a''$ of $\varrho$, we take the path $2a'$---$(a' + a'')$---$2a''$ in $G$. We denote the concatenation of these paths by $\varphi(\varrho)$. Intuitively, this operation scales $\varrho$ up by a factor of two. (Figure \ref{phi}.)

Observe that the length of $\varphi(\varrho)$ will be twice the length of $\varrho$, namely $2(k^2 - 1)$.

For each turn $a'ba''$ in $\varrho$, let us delete the subpath $(a' + b)$---$2b$---$(b + a'')$ from $\varphi(\varrho)$, and let us replace it with the edge $(a' + b)$---$(b + a'')$. We denote the resulting path in $G$ by $\psi(\varrho)$. Intuitively, this operation smooths down the sharp turns in $\varphi(\varrho)$. (Figure~\ref{psi}.)

Let $t$ be the number of turns in $\varrho$. Then the length of $\psi(\varrho)$ will be $2(k^2 - 1) - t$.

Finally, for each edge $a'a''$ of $\varrho$ which is in a free cycle, either we do nothing or, optionally, we choose one free cycle $a'a''b'b''$ which includes $a'a''$, we set $c = (a' + a'' + b' + b'')/4$, we delete the subpath $2a'$---$(a' + a'')$---$2a''$ from $\psi(\varrho)$, and we replace it with the subpath $2a'$---$2c$---$2a''$. We call a path in $G$ which can be obtained in this way a \emph{lift} of $\varrho$. Intuitively, this operation introduces some tiny aberrations in $\psi(\varrho)$. (Figure \ref{lift}.)

In particular, $\psi(\varrho)$ is also a lift of $\varrho$, namely the one in which we have selected the do-nothing option everywhere.

Observe that every lift of $\varrho$ is a snake path in $G$. Furthermore, all lifts of $\varrho$ are of the same length, namely $2(k^2 - 1) - t$.

We are ready to state and prove our structure theorem for the longest snake paths in~$G$.

\medskip

\begin{theorem} \label{king-odd-structure} Let $n$ be an odd positive integer with $n = 2k - 1$. Then every lift of a fewest-turn Hamiltonian path in the grid graph of size $k \times k$ is a longest snake path in the king graph of size $n \times n$. Conversely, every longest snake path in the king graph of size $n \times n$ can be obtained uniquely as a lift of some fewest-turn Hamiltonian path in the grid graph of size $k \times k$. \end{theorem}

\medskip

Theorem \ref{king-odd-structure} yields the following recipe for the generation of all longest snake paths in $G$: First, we generate all stamp-folding permutations of $k$ elements. Then we convert each stamp-folding permutation into four oriented fewest-turn Hamiltonian paths in $H$, two mostly-horizontal ones and two mostly-vertical ones. We forget about the orientations, and discard the duplicates. Finally, for each fewest-turn Hamiltonian path in $H$, we identify the corresponding free cycles, and we generate all of its lifts.

(We said in the introduction that each stamp-folding permutation yields two families of longest snake paths. Strictly speaking, the reality is that each permutation yields four families, and each quadruple of families is obtained in this way twice, out of two permutations $\sigma'$ and $\sigma''$ related by $\sigma'(i) + \sigma''(i) = k - 1$ for all $i$. However, it is straightforward to extract a one-to-two mapping from the two-to-four one.)

Before we go on to the proof of Theorem \ref{king-odd-structure}, let us briefly discuss the aberrations.

\medskip

\begin{proposition} \label{free} Let $\varrho$ be a fewest-turn Hamiltonian path in $H$ and let $f$ be the number of its corresponding free cycles. Then $\varrho$ yields exactly $2^f$ lifts. Furthermore, $f \le \max\{0, k - 5\}$, and for all $k$ this bound is attained by some $\varrho$. \end{proposition}

\medskip

\begin{proof} Suppose, for concreteness, that $\varrho$ is mostly-horizontal. Let $a'a''$ be an edge of $\varrho$ and let $a'a''b'b''$ be a free cycle which includes $a'a''$.

Suppose, for the sake of contradiction, that edge $a'a''$ is horizontal. Then, since $\varrho$ makes a turn at both of $b'$ and $b''$, the edges of $\varrho$ in the row of $b'$ and $b''$ cannot form one contiguous subpath of $\varrho$. This contradicts the fact that every row of the board contains exactly one horizontal segment of $\varrho$.

Thus edge $a'a''$ must be vertical.

Suppose now, for the sake of contradiction, that edge $a'a''$ is in a second free cycle $a'a''c'c''$. Since $\varrho$ must make a turn at both of $c'$ and $c''$ as well, it follows that the edges of $\varrho$ in the row of $b'$ and $c'$ cannot form one contiguous subpath of $\varrho$. From here, we get a contradiction as before.

Thus no edge of $\varrho$ can be in two distinct free cycles. Consequently, in the setting of Theorem \ref{king-odd-structure}, when we construct a lift of $\varrho$, we never have to choose between two free cycles which include the same edge $a'a''$.

We are left to show that $f \le \max\{0, k - 5\}$ and the bound is attained. We handle the cases when $k \le 6$ directly, and from now on we assume that $k \ge 7$.

Each free cycle contains two cells where $\varrho$ makes a turn. Conversely, each such cell is in at most one free cycle.

On the other hand, $\varrho$ makes $2k - 2$ turns altogether. However, a turn cell in an outermost column of $B$ cannot be in a free cycle. The lowermost turn cell and the topmost turn cell in a non-outermost column of $B$ cannot be in free cycles, either.

When at most one non-outermost column of $B$ contains turns of $\varrho$, since no turn cells outside of that column can be in free cycles, it follows that $2f \le k - 2$. Then $f \le k - 5$ by virtue of $k \ge 7$. Otherwise, when at least two non-outermost columns of $B$ contain turns of $\varrho$, it follows that at least eight turn cells are not in free cycles, and so $2f \le 2k - 10$.

The bound is attained, for example, when $\varrho$ corresponds to the stamp-folding permutation $0$, $2$, $3$, $4$, $\ldots$, $k - 1$, $1$. (The path in Figure \ref{rho} is of this form.) \end{proof}

\medskip

In the setting of the proof of Proposition \ref{free}, our observation that edge $a'a''$ must be vertical allows us to characterise the free cycles corresponding to $\varrho$ in terms of the underlying stamp-folding permutation $\sigma$, as follows:

Consider the ordered pairs $(\varepsilon, i)$ with $\varepsilon \in \{-1, 1\}$ and $0 \le i \le k - 2$ such that, in $\sigma$, both of $\sigma(i)$ and $\sigma(i + 1)$ lie between $\sigma(i) + \varepsilon(-1)^{\sigma(i)}$ and $\sigma(i + 1) + \varepsilon(-1)^{\sigma(i + 1)}$, and, additionally, there is some $j$ with $0 \le j \le k - 1$ such that all four of these lie between $j$ and $j + \varepsilon(-1)^j$.

Each such ordered pair yields a free cycle where edge $a'a''$ joins rows $i$ and $i + 1$. When $\varepsilon = 1$, cells $a'$ and $a''$ are in column $\omega_\text{Left}(\sigma(i)) - 1$ and cells $b'$ and $b''$ are on the right of them. Otherwise, when $\varepsilon = -1$, cells $a'$ and $a''$ are in column $k - \omega_\text{Right}(\sigma(i))$ and cells $b'$ and $b''$ are on their left. Furthermore, this accounts for all free cycles corresponding to~$\varrho$.

Let $\mathcal{N}_\text{Stamp}(k)$ be the number of stamp-folding permutations of $k$ elements and let $\mathcal{N}_\text{King}(n)$ be the number of longest snake paths in the king graph of size $n \times n$. Proposition~\ref{free} yields some loose bounds on $\mathcal{N}_\text{King}(n)$ in terms of $\mathcal{N}_\text{Stamp}(k)$.

\medskip

\begin{proposition} \label{log-bound} Let $n$ be an odd positive integer with $n = 2k - 1$. Then $2\mathcal{N}_\textnormal{Stamp}(k) < \mathcal{N}_\textnormal{King}(n) < 2^{k - 4}\mathcal{N}_\textnormal{Stamp}(k)$ for all $n$ with $n \ge 11$, and $\mathcal{N}_\textnormal{King}(n) = 2\mathcal{N}_\textnormal{Stamp}(k)$ when $3 \le n \le 9$. Thus, in particular, $\log \mathcal{N}_\textnormal{King}(n) = \Theta(n)$. \end{proposition}

\medskip

\begin{proof} The first part follows by Theorem \ref{king-odd-structure} and Proposition \ref{free}. (Strictly speaking, we should also note that not all fewest-turn Hamiltonian paths in $H$ attain the greatest number of free cycles when $n \ge 11$.) The second part is a corollary of the first part and the well-known asymptotic estimate $\log \mathcal{N}_\text{Stamp}(n) = \Theta(n)$. \cite{S} \end{proof}

\medskip

We continue with the proof of Theorem \ref{king-odd-structure}. Let $P$ be a longest snake path in $G$.

Here is a quick roadmap: Clearly, $P$ must attain exact equality in all inequalities from the proof of Theorem \ref{king-odd-path}. We examine all blocks of $G$ from this point of view, one by one, in a certain order, and we see that $P$ must satisfy certain purely local constraints. These constraints allow us to conclude that $P$ must be a lift.

The details, however, are somewhat technical.

We define an even cell $a$ of $A$ to be \emph{nice} (relative to $P$) when either $P$ visits $a$; or, else, $P$ does not visit $a$ but it traverses exactly one edge of $G$ between the four cells in the set $a + \{(1, 0), (0, 1), (-1, 0), (0, -1)\}$. (Note that, for some $a$, some of these cells might be outside of $A$.)

\medskip

\begin{lemma} \label{nice-i} Suppose that all even cells of $A$ are nice and $P$ does not visit any odd cells. Then there is a Hamiltonian path $\varrho$ in $H$ such that $P = \psi(\varrho)$. \end{lemma}

\medskip

\begin{proof} Let us call an edge of $G$ \emph{short} when the Euclidean distance between its endpoints is unity, and \emph{long} otherwise, when it is $\sqrt{2}$.

For each long edge $a'a''$ of $P$, we do the following: Since $P$ does not visit any odd cells, there is a unique even cell $b$ such that both of $a'b$ and $ba''$ are edges of $G$. We delete edge $a'a''$ from $P$, and we replace it with these two edges.

Since all even cells of $A$ are nice and $P$ is a snake path, the result will be a path in $G$ which contains only short edges, which visits all even cells of $A$, and which does not visit any odd cells. Denote this path by $Q$.

Let $c$ be an endpoint of $Q$. Suppose, for the sake of contradiction, that $c$ is not an even cell. Since $Q$ does not visit any odd cells, there are exactly two even cells $d'$ and $d''$ of $A$ adjacent to $c$ in $G$. Let $cd'$ be the unique edge of $Q$ incident with $c$. Observe that $c$ is also an endpoint of $P$. But then $d''$ cannot be nice because $P$ is a snake path, and we arrive at a contradiction.

So, to our previous observations about $Q$, we can add the fact that both of its endpoints are even cells. Consequently, $Q$ is of the form $Q = \varphi(\rho)$ for some Hamiltonian path $\varrho$ in $H$, and $P = \psi(\varrho)$, as needed. \end{proof}

\medskip

\begin{figure}[ht] \centering \includegraphics{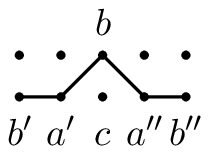} \caption{} \label{rectify} \end{figure}

We define a \emph{rectifiable aberration} in $P$ to be a subpath of $P$ of the form $b'a'ba''b''$ such that $a'$ and $a''$ are two even cells in the same row or column; the cell $c = (a' + a'')/2$ is a common neighbour of $a'$ and $a''$ in $G$; the cells $b'$ and $b''$ satisfy $a' = (b' + c)/2$ and $a'' = (c + b'')/2$; and $b \neq c$. (Figure \ref{rectify}.)

To \emph{rectify} a rectifiable aberration, we delete the subpath $a'ba''$ from $P$, and we replace it with the subpath $a'ca''$. The result will be a new snake path in $G$ of the same length as $P$. This follows because $P$ being a snake path implies that the only neighbours of $c$ in $G$ which $P$ visits are $a'$, $b$, and $a''$.

\medskip

\begin{lemma} \label{nice-ii} Suppose that $P$ does not contain any rectifiable aberrations. Then all even cells of $A$ are nice and $P$ does not visit any odd cells. \end{lemma}

\medskip

\begin{proof} Since $P$ is a longest snake path in $G$, it must attain exact equality in all inequalities from the proof of Theorem \ref{king-odd-path}. Thus:

(i) A cell of $A$ of the form $(z, z)$ with $z$ odd cannot be in $P$ because it is an odd cell which is in more than two blocks. Same goes for the images of these cells under the symmetries of $A$;

(ii) An edge of $G$ of the form $(z, z + 1)$---$(z + 1, z)$ with $z$ odd and $z \le k - 2$ cannot be in $P$ because it is a regular edge which is in more than one block. Same goes for the images of these edges under the symmetries of $A$; and

(iii) Every block $\nosovka$ must attain exact equality in Lemma \ref{blocks}, so that $w(\nosovka) = 1$ when $\nosovka$ is a little block and $w(\nosovka) = 2$ otherwise, when $\nosovka$ is a large block.

For all odd positive integers $s$ with $1 \le s \le n$, we write $A_s$ for the concentric subboard of $A$ of size $s \times s$ given by $A_s = [\lfloor n/2 \rfloor - \lfloor s/2 \rfloor; \lfloor n/2 \rfloor + \lfloor s/2 \rfloor]^2$.

We will show by induction on $s$ that all even cells of $A_s$ are nice and all odd cells of $A_s$ are outside of $P$.

Our base case is $s = 1$. Then $A_1$ consists of a single cell, namely $(k - 1, k - 1)$. Denote this cell by $o$.

When $k$ is even, $o$ is an odd cell, and $o \not \in P$ by (i).

\begin{figure}[ht] \centering \includegraphics{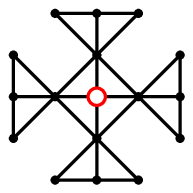} \caption{} \label{clover} \end{figure}

When $k$ is odd, $o$ is an even cell. (Figure \ref{clover}.) Suppose, for the sake of contradiction, that $k \ge 3$ and $o$ is not nice. Then $o \not \in P$. By (i), $o + (1, 1) \not \in P$, and similarly for the images of this cell under the symmetries of $A$. By (ii), $(o + (1, 0))$---$(o + (0, 1)) \not \in P$, and similarly for the images of this edge under the symmetries of $A$. Thus (iii) for $\nosovka(o + (0, -2))$ implies $o + (-1, -2) \in P$, $o + (1, -2) \in P$, and either $o + (0, -1) \in P$ or $o + (0, -2) \in P$; and similarly for the images of this block under the symmetries of $P$. However, the twelve cells we just concluded must be in $P$ are the vertices of a cycle in $G$, and we arrive at a contradiction.

This settles the base case.

For the induction step, let $s \ge 3$ and suppose that we have already established the desired result for $A_{s - 2}$. Let $a = (x, y)$ be an arbitrary cell of $A_s \setminus A_{s - 2}$. By symmetry, we can assume without loss of generality that $\lfloor n/2 \rfloor - \lfloor s/2 \rfloor \le x \le \lfloor n/2 \rfloor$ and $y = \lfloor n/2 \rfloor - \lfloor s/2 \rfloor$.

Figures \ref{step-cases-ab} and \ref{step-cases-cd} show some of the cells, edges, and blocks relevant to our reasoning. Cell $a$ is highlighted in all of them. Note that some blocks which are shown as large in the figures might be little ones in reality. In all such cases, we emphasise this possibility in the text.

Suppose, for the sake of contradiction, that $a$ is an even cell but that it is not nice.

\begin{figure}[t] \null \hfill \begin{subfigure}{45pt} \centering \includegraphics{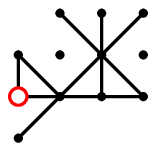} \caption{} \label{step-case-a} \end{subfigure} \hfill \begin{subfigure}{55pt} \centering \includegraphics{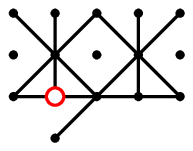} \caption{} \label{step-case-b} \end{subfigure} \hfill \null \caption{} \label{step-cases-ab} \end{figure}

Then $a \not \in P$. Furthermore, the odd cells of $\nosovka(a)$ are not in $P$, either, by the induction hypothesis.

\smallskip

\emph{Case 1}. $x = \lfloor n/2 \rfloor - \lfloor s/2 \rfloor$ and $\nosovka(a)$ is a little block. (Figure \ref{step-case-a}.)

Then (iii) for $\nosovka(a)$ implies $(a + (1, 0))$---$(a + (0, 1)) \in P$. Since $a$ is not nice, at least one more edge of $G$ between the four cells in the set $a + \{(1, 0), (0, 1), (-1, 0), (0, -1)\}$ must be in $P$. It cannot be $(a + (-1, 0))$---$(a + (0, -1))$ because then $P$ would contain the vertices of a cycle in $G$. The other two subcases are symmetric with respect to the line of unit slope through $a$, and we assume that $(a + (0, -1))$---$(a + (1, 0)) \in P$.

It follows that $a + (2, 0) \not \in P$ and $a + (2, 1) \not \in P$. Thus no edges of $\nosovka(a + (2, 0))$ are in $P$. Since also the odd cells of this block are outside of $P$ by the induction hypothesis, we arrive at a contradiction with (iii). (The conclusion holds regardless of whether the block is a little one or a large one.)

\smallskip

\emph{Case 2}. $\lfloor n/2 \rfloor - \lfloor s/2 \rfloor < x \le \lfloor n/2 \rfloor$ and $\nosovka(a)$ is a large block. (Figure \ref{step-case-b}.)

By (iii) for $\nosovka(a)$ and the induction hypothesis for $a + (0, 2)$, exactly one of the two edges $(a + (1, 0))$---$(a + (0, 1))$ and $(a + (0, 1))$---$(a + (-1, 0))$ is in $P$. The two subcases are analogous, and we assume that the former edge is in $P$ while the latter one is not. As in Case 1, at least one more edge of $G$ between the four cells in the set $a + \{(1, 0), (0, 1), (-1, 0), (0, -1)\}$ must be in $P$, and it cannot be $(a + (-1, 0))$---$(a + (0, -1))$ because then $P$ would contain the vertices of a cycle in $G$. Thus $(a + (0, -1))$---$(a + (1, 0)) \in P$.

From here, we arrive at the exact same contradiction as in Case 1. (Once again, regardless of the type of the block $\nosovka(a + (2, 0))$.)

\smallskip

We have established that if $a$ is an even cell, then it is nice.

For the second half of the induction step, suppose, for the sake of contradiction, that $a$ is an odd cell with $a \in P$. By (i), it follows that $\lfloor n/2 \rfloor - \lfloor s/2 \rfloor + 2 \le x \le \lfloor n/2 \rfloor$.

\smallskip

\emph{Case 1}. $a + (-2, 0) \in P$. (Figure \ref{step-case-c}.)

By (i), $x \ge \lfloor n/2 \rfloor - \lfloor s/2 \rfloor + 4$. Thus $\nosovka(a + (-1, 1))$ is a large block.

Note that $a \in P$ implies $(a + (-1, 1))$---$(a + (0, 1)) \not \in P$ and $a + (-2, 0) \in P$ implies $(a + (-1, 1))$---$(a + (-2, 1)) \not \in P$. Furthermore, by the induction hypothesis, the odd cells of $\nosovka(a + (-1, 1))$ are not in $P$ and cell $a + (-1, 3)$ is nice. Then (iii) for $\nosovka(a + (-1, 1))$ implies that exactly one of the two edges $(a + (0, 1))$---$(a + (-1, 2))$ and $(a + (-1, 2))$---$(a + (-2, 1))$ is in $P$. (Here, we take into account the fact that $a \in P$ and $a + (-2, 0) \in P$ together imply $(a + (-1, 1))$---$(a + (-1, 2)) \not \in P$.) The two subcases are analogous, and we assume that the former edge is in $P$ while the latter one is not.

\begin{figure}[t] \null \hfill \begin{subfigure}{55pt} \centering \includegraphics{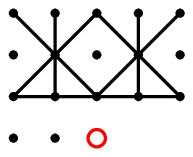} \caption{} \label{step-case-c} \end{subfigure} \hfill \begin{subfigure}{100pt} \centering \includegraphics{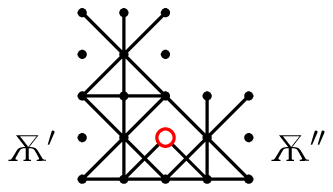} \caption{} \label{step-case-d} \end{subfigure} \hfill \null \caption{} \label{step-cases-cd} \end{figure}

It follows that $a$---$(a + (0, 1))$---$(a + (-1, 2)) \subseteq P$, $a + (1, 1) \not \in P$, and $a + (1, 2) \not \in P$. Thus no edges of $\nosovka(a + (1, 1))$ are in $P$. Since also the odd cells of this block are outside of $P$ by the induction hypothesis, we arrive at a contradiction with (iii). (The block $\nosovka(a + (1, 1))$ will always be a large one because of our symmetry-breaking assumption that $x \le \lfloor n/2 \rfloor$. However, in the analogous subcase when $(a + (0, 1))$---$(a + (-1, 2)) \not \in P$ and $(a + (-1, 2))$---$(a + (-2, 1)) \in P$, the contradiction occurs at block $\nosovka(a + (-3, 1))$ which could happen to be a little one.)

\smallskip

\emph{Case 2}. $a + (2, 0) \in P$. This case is analogous to Case 1.

\smallskip

\emph{Case 3}. $a + (-2, 0) \not \in P$ and $a + (2, 0) \not \in P$. (Figure \ref{step-case-d}.)

Observe that both of $\nosovka' = \nosovka(a + (-1, -1))$ and $\nosovka'' = a + (1, -1)$ are large blocks.

Since $a \in P$, all edges of $\nosovka'$ other than $(a + (-1, -1))$---$(a + (-2, -1))$, $(a + (-2, -1))$---$(a + (-1, 0))$, and $(a + (-1, 0))$---$(a + (-2, 1))$ are not in $P$. Then, in light of $a \in P$ and $a + (-2, 0) \not \in P$, (iii) for $\nosovka'$ implies that exactly one of these edges is in $P$. In particular, exactly one cell $b'$ out of the pair $a + \{(-1, -1), (-1, 0)\}$ is in $P$. Similar reasoning applies to $\nosovka''$, and we define $b''$ analogously.

Suppose, for the sake of contradiction, that $b' = a + (-1, 0)$. Then $ab' \subseteq P$ implies $a + (0, 1) \not \in P$ and $a + (-1, 1) \not \in P$. Furthermore, the odd cells of $\nosovka(a + (-1, 1))$ are outside of $P$ by the induction hypothesis. When this block is a little one, we arrive at a contradiction with (iii) immediately. Otherwise, when it is a large one, we arrive at a contradiction with (iii) anyway once we take into account the fact that, by the induction hypothesis, both of $a + (-1, 1)$ and $a + (-1, 3)$ are nice.

Consequently, $b' = a + (-1, -1)$. Similarly, $b'' = a + (1, -1)$. By (iii) for $\nosovka'$ and $\nosovka''$, it follows that also $(a + (-1, -1))$---$(a + (-2, -1)) \in P$ and $(a + (1, -1))$---$(a + (2, -1)) \in P$. However, then $(a + (-2, -1))$---$(a + (-1, -1))$---$a$---$(a + (1, -1))$---$(a + (2, -1))$ becomes a rectifiable aberration in $P$, and we arrive at a contradiction. (In fact, this is the only place in the proof where we use the constraint that $P$ does not contain any rectifiable aberrations.)

\smallskip

We have established that if $a$ is an odd cell, then it cannot be in $P$. The induction step is complete. \end{proof}

\medskip

We are ready to tackle Theorem \ref{king-odd-structure}.

\medskip

\begin{proof*}{Proof of Theorem \ref{king-odd-structure}} For a start, let us rectify all rectifiable aberrations in $P$ one by one. The result will be a snake path $Q$ in $G$ of the same length as $P$ and without any rectifiable aberrations.

By Lemma \ref{nice-ii}, we get that all even cells are nice relative to $Q$ and $Q$ does not visit any odd cells.

By Lemma \ref{nice-i}, it follows that there is a Hamiltonian path $\varrho$ in $H$ such that $Q = \psi(\varrho)$.

Let $t$ be the number of turns in $\varrho$. Then both of $P$ and $Q$ are of length $2(k^2 - 1) - t$.

Since the greatest length of a snake path in $G$ is $(n^2 - 1)/2$ by Theorem \ref{king-odd-path}, we conclude that $t = 2k - 2$, and so $\varrho$ is in fact a fewest-turn Hamiltonian path in $H$.

Finally, in order to transform $Q$ back into $P$, we must restore the rectifiable aberrations which we removed in the beginning. However, it is straightforward to check that the spots in $Q$ where we can introduce a rectifiable aberration are exactly the ones associated with the free cycles corresponding to $\varrho$. Therefore, $P$ is a lift of $\varrho$, as needed.

The reasoning in the last few paragraphs shows also the converse: That every lift of a fewest-turn Hamiltonian path in $H$ is a longest snake path in $G$. \end{proof*}

%% file: snake-06-knight.tex
\section{Knight Graphs} \label{knight}

In this section, we prove Theorem \ref{knight-all}. Let $m$ and $n$ be positive integers, let $A$ be the standard board of size $m \times n$, and let $G$ be the knight graph on $A$.

We begin with the upper bound.

One natural approach would be as follows: First we find some finite knight graph $H$ with pseudosnake density $1/2$. Then we sum over all translation copies of $H$ contained within our board.

The author was not able to implement this idea in its purest form. Below, we present a slightly more complicated argument which relies on a weighted knight graph instead.

We define a \emph{weighted graph} $\Gamma$ to consist of a simple graph $H$ and a weighting function $w$ which assigns a nonnegative real weight to each vertex of $H$. When $H$ is finite, we denote the total weight of all of its vertices by $w(\Gamma)$, and we define the pseudosnake density of $\Gamma$ to be the ratio of the greatest total weight of the vertices in a pseudosnake of $H$ to $w(\Gamma)$.

\medskip

\begin{lemma} \label{capacity} Suppose that there is a weighted knight graph $\Gamma$ with pseudosnake density~$\tau$. Then the number of vertices in a pseudosnake of $G$ cannot exceed $\tau mn + \mathcal{O}(m + n)$. \end{lemma}

\medskip

Here and in the proof, the implicit constants in the $\mathcal{O}$-terms depend on $\Gamma$.

\medskip

\begin{proof} Let $P$ be a pseudosnake of $G$.

Consider all translation copies of $\Gamma$ that fit within $A$. There are $mn + \mathcal{O}(m + n)$ of them. For each such copy, the total weight of all cells of $P$ within it cannot exceed~$\tau w(\Gamma)$.

Conversely, it is true of all but $\mathcal{O}(m + n)$ cells $a$ of $P$ that $a$ is sufficiently far away from the boundary of $A$ for every translation copy of $\Gamma$ which contains $a$ to fit within $A$. For each such cell of $P$, the sum of its weights over all translation copies of $\Gamma$ which contain it will be $w(\Gamma)$. \end{proof}

\medskip

\begin{proof*}{Proof of the upper bound for Theorem \ref{knight-all}} By Lemma \ref{capacity}, it suffices to exhibit one concrete weighted knight graph with pseudosnake density $1/2$. We claim that the one in Figure \ref{pancake} works. (The figure shows all cells of the graph together with the weights assigned to them.) It has $68$ cells of total weight $192$.

\begin{figure}[ht] \centering \includegraphics{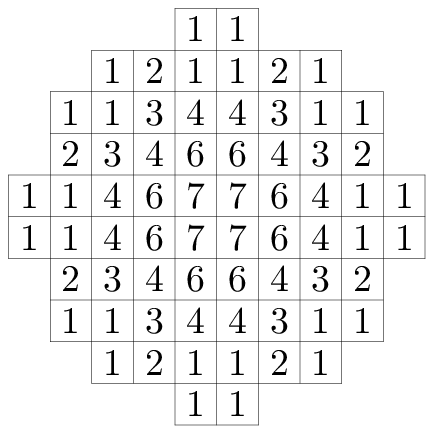} \caption{} \label{pancake} \end{figure}

Our claim would likely be extremely difficult to check by hand. However, it is straightforward to check with the help of a standard constraint satisfaction solver. The author has done this twice, using two different constraint satisfaction frameworks: the Copris package for the Scala programming language and the OR-Tools package for the Python programming language. \end{proof*}

\medskip

One might wonder how the weighted knight graph in Figure \ref{pancake} was found.

\begin{figure}[ht] \null \hfill \begin{minipage}{80pt} \centering \includegraphics{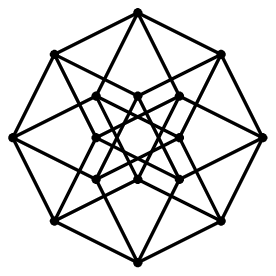} \caption{} \label{hypercube} \end{minipage} \hfill \begin{minipage}{45pt} \centering \includegraphics{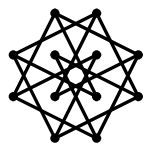} \caption{} \label{ball} \end{minipage} \hfill \null \end{figure}

For a start, let $S_\text{Tess}$ be the set of all $16$ cells of the form $\varepsilon_1(2, 1) + \varepsilon_2(1, 2) + \varepsilon_3(-1, 2) + \varepsilon_4(-2, 1)$, where $\varepsilon_i \in \{0, 1\}$ for all $i$. Then the knight graph $G_\text{Tess}$ on $S_\text{Tess}$ is isomorphic to the tesseract graph. (Figure \ref{hypercube}.)

It is not too difficult to check by hand that the pseudosnake density of the tesseract graph is $9/16$, and that it is attained by an essentially unique pseudosnake. Since $9/16$ is very close to $1/2$ from above, we see that $G_\text{Tess}$ works almost, but not quite.

We can attempt to fix this by taking the union of several overlapping copies of $G_\text{Tess}$. Since $G_\text{Tess}$ itself just barely manages to push through pseudosnake density $1/2$, we can hope that the interference between its copies will prevent too many of them from doing the same.

We formalise this notion as follows: Let $S'$ and $S''$ be two nonempty finite sets of cells. We define their \emph{sum}, denoted $S' + S''$, to be the multiset of cells which consists of all cells of the form $a' + a''$ with $a' \in S'$ and $a'' \in S''$, and where the multiplicity of each cell is the number of ways that it can be expressed in this form.

Then, given a multiset of cells $\mathcal{S}$, we define $\knightgraph{\mathcal{S}}$ to be the weighted knight graph on the cells of $\mathcal{S}$ where the weight of each cell is its multiplicity in $\mathcal{S}$.

For each nonempty finite set of cells $S$, we can think of the weighted knight graph $\knightgraph{S + S_\text{Tess}}$ as constructed out of several overlapping copies of $G_\text{Tess}$. We experiment with different $S$, and eventually we strike gold with $S_\text{Dia} = [0; 3]^2 \setminus \{0, 3\}^2$. This is the Aztec diamond of order two. (Figure \ref{ball}.)

We go on to the lower bound.

One natural approach would be as follows: First we find a doubly periodic pseudosnake $P_\infty$ in $\knightgraph{\mathbb{Z}^2}$ with density $1/2$, where furthermore every cell is of degree exactly two and there are no finite cycles. (Here, ``doubly periodic'' means that there are two linearly independent two-dimensional vectors $u$ and $v$ with $a \in P_\infty \Leftrightarrow a + u \in P_\infty \Leftrightarrow a + v \in P_\infty$ for all cells $a$ of $\mathbb{Z}^2$.)

Then, given a board $A$, we take the restriction $P^\star$ of $P_\infty$ to $A$. Because of the structure of $P_\infty$, this restriction will be the disjoint union of several paths. We make some modifications near the boundary of $A$, deleting some cells from $P^\star$ and replacing them with new ones, so as to stitch all of these paths together into a single snake path or cycle. Since $P_\infty$ is doubly periodic with density $1/2$, originally $P^\star$ will contain $mn/2 + \mathcal{O}(m + n)$ cells. Finally, keeping our modifications close to the boundary of $A$ ensures that they cost us $\mathcal{O}(m + n)$ of these cells altogether.

Finding a suitable $P_\infty$ is straightforward enough. For example, the set of all cells $(x, y)$ with $x \bmod 4 \in \{0, 1\}$ works. It consists of vertical strips of width two spaced two units apart.

However, the second part of our plan runs into significant difficulties. Thus the construction we present below is somewhat more complicated. We divide $A$ into four large regions; we fill up different regions using different pseudosnakes $P_\infty$; and we make stitching-together modifications not only near the boundary of $A$, but also near the boundaries between regions.

We continue with the details.

We define a \emph{twine} to be a board of height two. When $s \ge 2$, the knight graph on a twine with width $s$ is the disjoint union of four paths, two spanning $\lfloor s/2 \rfloor$ cells each and two spanning $\lceil s/2 \rceil$ cells each.

\begin{figure}[ht] \centering \null \hfill \begin{minipage}{90pt} \centering \includegraphics{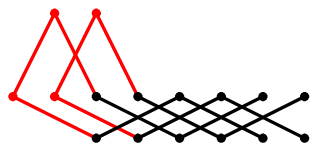} \caption{} \label{tie} \end{minipage} \hfill \begin{minipage}{115pt} \centering \includegraphics{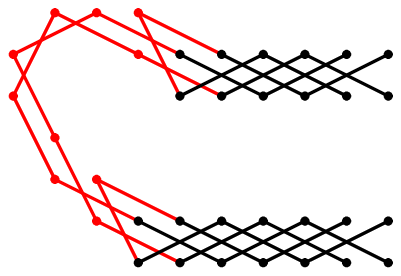} \caption{} \label{splice} \end{minipage} \hfill \null \end{figure}

Consider a twine $E$ with lower left corner cell $a$. To \emph{tie off} $E$ on the left, we add to it the four cells in the set $a + \{(-2, 1), (-1, 1), (-1, 3), (0, 3)\}$. (Figure \ref{tie}.) Similarly, to tie off on the right a twine $E$ with lower right corner $a$, we add to it the reflections of the four cells above with respect to the vertical line through $a$.

Consider, now, two twines $E$ and $F$ with lower left corners $a$ and $b$ satisfying $a + (1, 4) = b$. To \emph{splice together} $E$ and $F$ on the left, we add to them the ten cells in the set $a + \{(-3, 4), (-3, 5), (-2, 2), (-2, 3), (-2, 6), (-1, 1), (-1, 2), (-1, 6), (0, 5), (0, 6)\}$. (Figure \ref{splice}.) Similarly, to splice together on the right two twines $E$ and $F$ whose lower right corners $a$ and $b$ satisfy $a + (-1, 4) = b$, we add to them the reflections of the ten cells above with respect to the vertical line through $a$.

Let $k$ be a positive integer and let $I = [x'; x'']$ be an integer interval with $|I| \ge 8k - 5$. For each $i$ with $0 \le i \le k - 1$, if $i$ is even then construct the twine $E_i = [x' + 4i; x'' - 4i] \times [4i; 4i + 1]$, and if $i$ is odd then construct the twine $E_i = [x' + 4i + 3; x'' - 4i + 3] \times [4i; 4i + 1]$. Tie off $E_0$ on the left; for all $i$ with $0 \le i \le k - 2$, splice together $E_i$ and $E_{i + 1}$ on the right if $i$ is even, and on the left if $i$ is odd; and, finally, if $k$ is even then tie off $E_{k - 1}$ on the left, and if $k$ is odd then tie it off on the right.

\begin{figure}[ht] \centering \includegraphics{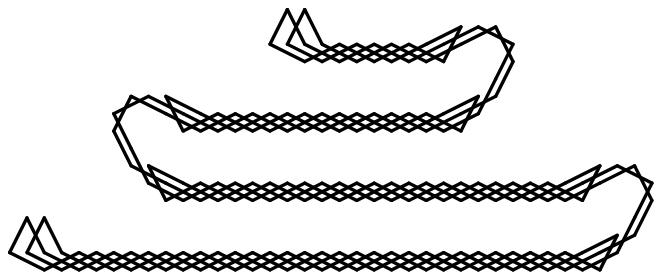} \caption{} \label{skein} \end{figure}

We denote the resulting set of cells by $U(k, I)$. For example, Figure \ref{skein} shows the knight graph on $U(4, [0; 32])$.

\medskip

\begin{lemma} \label{skein-cycle} Suppose that $k \ge 2$ and $|I|$ is odd. Then the knight graph on $U(k, I)$ is a cycle. \end{lemma}

\medskip

\begin{proof} Denote $H = \knightgraph{U(k, I)}$. It is straightforward to check that all cells of $H$ are of degree two. We are left to verify that $H$ is connected.

Suppose, for the sake of clarity, that $k$ is odd. The opposite case, when it is even, is similar.

Let $a_i$ and $b_i$ be the lower left and lower right corner cells of $E_i$. Let also $a'_i$ and $a''_i$ be the two cells of $E_i$ adjacent by side to $a_i$, and define $b'_i$ and $b''_i$ similarly for $b_i$.

It is straightforward to verify that, since $|I|$ is odd: (a) A path in $H$ connects $b'_0$ and $b''_0$ and covers $E_0$ together with the cells which tie it off; (b) For each $i$ with $1 \le i \le k - 2$, two paths in $H$ connect the pairs $\{a'_i, a''_i\}$ and $\{b'_i, b''_i\}$ and cover $E_i$ together with two cells in the adjacent splices; and (c) A path in $H$ connects $a'_{k - 1}$ and $a''_{k - 1}$ and covers $E_{k - 1}$ together with the cells which tie it off.

Finally, the remaining cells of the splices of $U(k, I)$ form additional knight paths in $H$ which connect the pairs $\{b'_i, b''_i\}$ and $\{b'_{i + 1}, b''_{i + 1}\}$ for all even $i$ with $0 \le i \le k - 3$ as well as the pairs $\{a'_i, a''_i\}$ and $\{a'_{i + 1}, a''_{i + 1}\}$ for all odd $i$ with $1 \le i \le k - 2$. \end{proof}

\medskip

\begin{proof*}{Proof of the lower bound for Theorem \ref{knight-all}} We construct a large snake cycle in $G$, and for a large snake path we can simply delete one cell from that cycle.

Suppose, without loss of generality, that $m \le n$. Since we are already willing to accept a tolerance of $\mathcal{O}(m + n)$, we can safely assume that $m = 8k + 14$ for some positive integer $k$ with $k \ge 3$, and that $n$ is even.

Construct $U_\text{I} = U(k, [8; n - 12])$. Let also $V_\text{I}$ be the reflection of $U(k, [8; m - 12])$ with respect to the line $x = y$. Lastly, let $U_\text{II}$ and $V_\text{II}$ be symmetric to $U_\text{I}$ and $V_\text{I}$ with respect to the center of $A$.

The knight graph on $U_\text{I} \cup V_\text{I} \cup U_\text{II} \cup V_\text{II}$ is the disjoint union of four cycles. We proceed to stitch these four cycles together into a single longer cycle.

\begin{figure}[ht] \centering \includegraphics{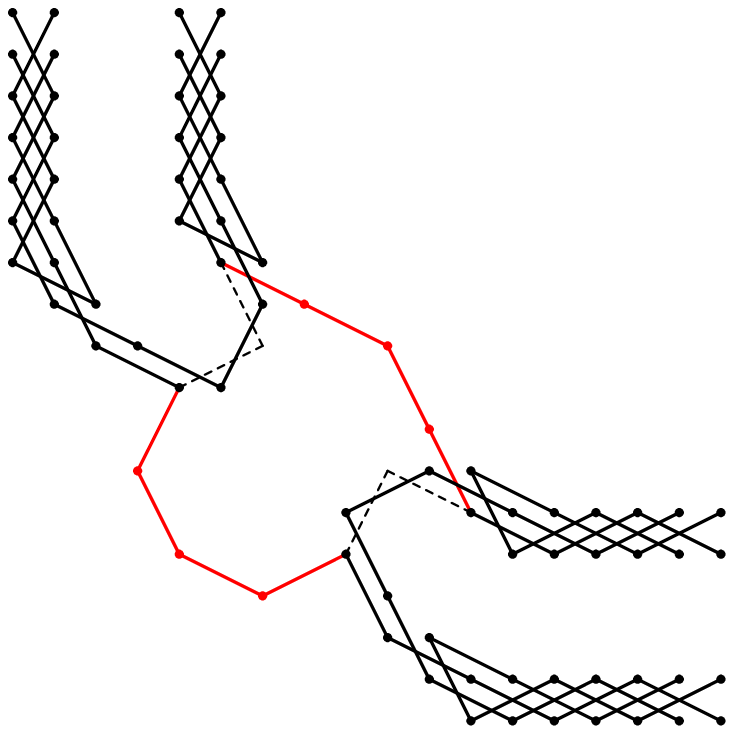} \caption{} \label{stitch} \end{figure}

Define $S_\text{Del} = \{(6, 9), (9, 6)\}$ and $S_\text{Add} = \{(3, 6), (4, 4), (6, 3), (7, 10), (9, 9), (10, 7)\}$. Delete the two cells of $(4, 4) + S_\text{Del}$ from $U_\text{I}$ and $V_\text{I}$, and replace them with the six cells of $(4, 4) + S_\text{Add}$. This stitches together the cycles of $U_\text{I}$ and $V_\text{I}$. (Figure \ref{stitch}.)

We carry out two more such modifications. For one of them, we reflect the sets $S_\text{Del}$ and $S_\text{Add}$ with respect to the vertical axis of symmetry of the board, we delete the two cells in the image of $S_\text{Del}$ from $U_\text{I}$ and $V_\text{II}$, and we replace them with the six cells in the image of $S_\text{Add}$. This stitches together the cycles of $U_\text{I}$ and $V_\text{II}$. For the other one, we proceed similarly, except that the reflections are done with respect to the horizontal axis of symmetry of the board. This stitches together the cycles of $V_\text{I}$ and $U_\text{II}$.

\begin{figure}[ht] \centering \includegraphics{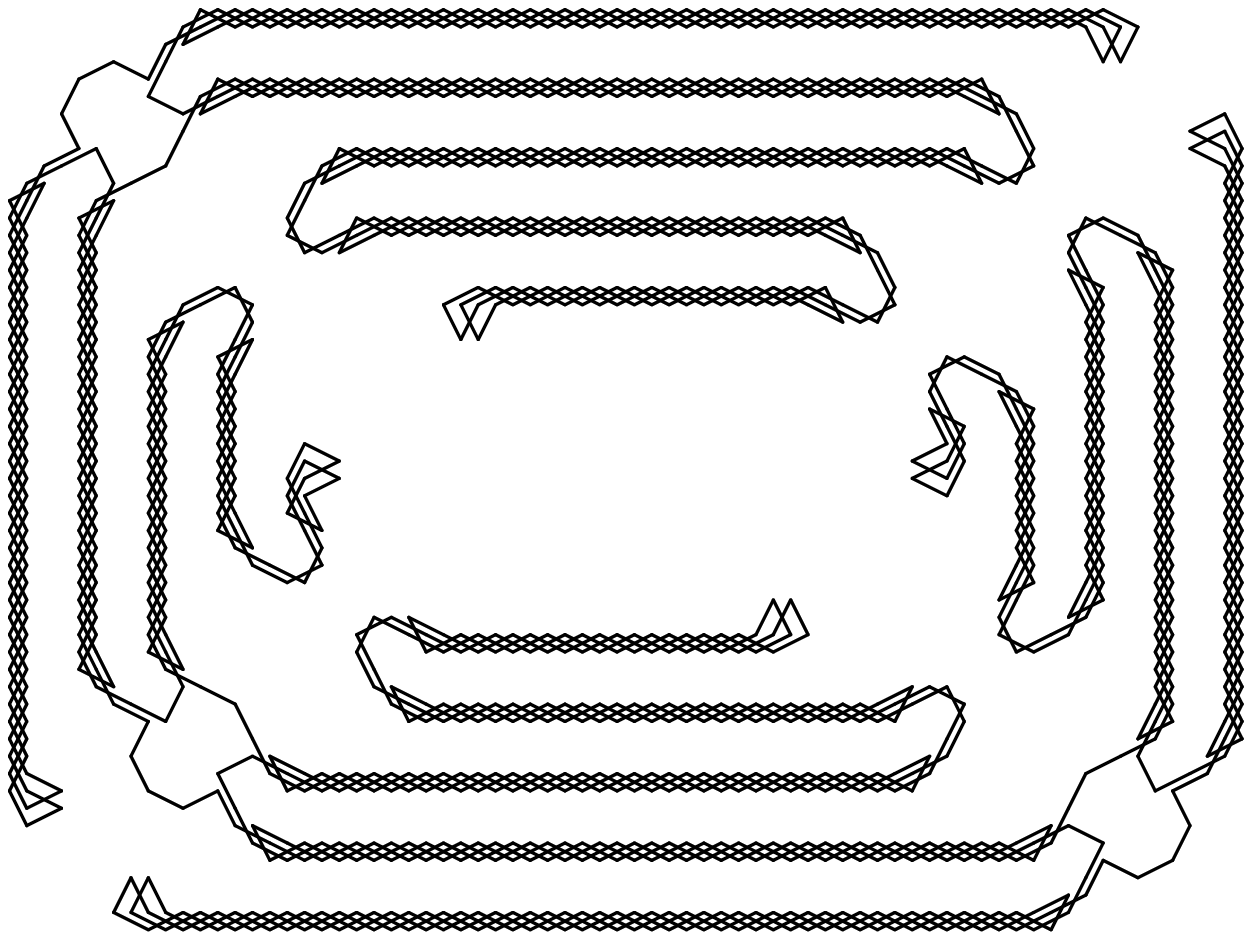} \caption{} \label{knight-large} \end{figure}

Let $W$ be the final set of cells obtained in this way. For example, Figure \ref{knight-large} shows $W$ in the case when $k = 5$, $m = 54$, and $n = 72$.

Observe that the density of $W$ within $A$ is $1/2$ everywhere except within five strips of bounded width. (Four of these strips surround portions of the interior angle bisectors at the four corners of $A$, and the fifth one surrounds a portion of the horizontal axis of symmetry of $A$. The corresponding mostly hollow areas are clearly visible in Figure \ref{knight-large}.) Consequently, the number of cells in $W$ is $mn/2 + \mathcal{O}(m + n)$.

On the other hand, $W$ is the vertex set of a snake cycle in $G$. \end{proof*}

%% file: snake-07-further-i.tex
\section{Further Work on King Graphs} \label{further-i}

In this section, we collect some additional results and open problems on king graphs.

We saw that the behaviour of the longest snake paths in $\kinggraph{n \times n}$ depends on the parity of $n$. For cycles, it appears that there are four classes instead, depending on the value of $n \bmod 4$. The techniques we developed for paths quickly resolve two of them.

\medskip

\begin{theorem} \label{king-0mod4} Let $n$ be a positive integer with $n \equiv 0 \pmod 4$ and $n \ge 8$. Then the greatest length of a snake cycle in the king graph of size $n \times n$ is $n^2/2 - 1$. Furthermore, for all such $n$, there are exactly $48$ snake cycles which attain this greatest length. These cycles are all asymmetric, and so six of them are essentially distinct. \end{theorem}

\medskip

For completeness, there is a unique snake cycle of the greatest length $8$ when $n = 4$.

Each one of the cycles of Theorem \ref{king-0mod4} is shaped like a double spiral.

It is curious that the number of longest snake cycles freezes in this way, and the cycles themselves crystallise into a single inflexible structure. A similar phenomenon occurs in the setting of Theorem \ref{near-hamiltonian}.

\medskip

\begin{proof} Let $n = 2k$ and define $A$, $B$, $G$, $H$, and $\Phi$ as in Section \ref{king-even}. Let also $C$ be a snake cycle in $G$ of length at least $n^2/2 - 1$.

As in Section \ref{king-even}, for each cell $b$ of $B$ at most two cells of $\Phi(b)$ are in $C$. Let us call $b$ \emph{deficient} when this bound is not attained. Thus there is at most one deficient cell.

Observe that, if an edge of $C$ joins one cell of $\Phi(b')$ and one cell of $\Phi(b'')$, with $b' \neq b''$, then $b'b''$ must still be an edge of $H$. Otherwise, assuming for concreteness that $b' + (1, 1) = b''$, by an argument similar to the one in Section \ref{king-even} we see that at least two of the four cells in the set $b' + [0; 1]^2$ must be deficient, a contradiction.

This allows us to define the Hamiltonian cycle $\varrho$ in $H$ relative to $C$ in the same manner as in Section \ref{king-even}.

Define also the cycle $E$ in $H$ as in the proof of Lemma \ref{spiral}. Since $E$ is not a Hamiltonian cycle of $H$ when $k \ge 4$, at least one edge of $E$ must be outside of $\varrho$.

On the other hand, observe that Lemma \ref{turn} now admits a unique exception: When the turn occurs at a deficient cell. Consequently, every edge of $E$ outside of $\varrho$ must possess a deficient endpoint.

Define $A^\star$, $B^\star$, $G^\star$, and $H^\star$ as in Section \ref{king-even}. Let also $C^\star$ and $\varrho^\star$ be the restrictions of $C$ and $\varrho$ to $A^\star$ and $B^\star$, respectively. We obtain that: (a) There is exactly one deficient cell; (b) There is exactly one edge $\beta'\beta''$ of $E$ outside of $\varrho$; (c) $\varrho^\star$ is a Hamiltonian path in $H^\star$ whose endpoints are the two neighbours of $\beta'$ and $\beta''$ in $B^\star$; and (d) $C^\star$ is a snake path in $G^\star$ of the greatest possible length whose associated Hamiltonian path in $H^\star$ is~$\varrho^\star$.

However, the proof of Theorem \ref{king-even-path} gives us a complete description of all Hamiltonian paths in $H^\star$ associated with a longest snake path in $G^\star$. Since the two endpoints of $\varrho^\star$ are neighbours in $H^\star$, we conclude that when $n \ge 12$ a symmetry of $B^\star$ must map $\varrho^\star$ onto the unique regular Hamiltonian path in $H^\star$ of type $\operatorname{II}(0)$, as defined in Section \ref{king-even}. The rest is straightforward. \end{proof}

\medskip

The other class we can tackle without too much extra effort is $n \equiv 3 \pmod 4$. First, though, we need to sort through some preliminaries.

Let $s$ be an even positive integer. Given a permutation $\sigma$ of $[0; s - 1]$, consider a closed curve in the plane defined in the same way as the curve $\kappa$ in Section \ref{king-odd-ii}, except that one additional arc on the right of the coordinate axis $Oy$ joins points $(0, \sigma^{-1}(s - 1))$ and $(0, \sigma^{-1}(0))$. A curve of this form which does not intersect itself is known as a \emph{closed meander}.

Consider a Hamiltonian cycle in the grid graph $\Gamma$ of size $s \times s$. The smallest number of turns that such a cycle can make is $2s$. \cite{J} Furthermore, the closed meanders with $s$ arcs and the fewest-turn Hamiltonian cycles in $\Gamma$ are related in a way analogous to the relation between stamp-folding permutations and fewest-turn Hamiltonian paths.

Observe, lastly, that our definition of a lift in Section \ref{king-odd-ii} works just as well with cycles instead of paths.

\medskip

\begin{theorem} \label{king-3mod4} Let $n$ be a positive integer with $n \equiv 3 \pmod 4$ and $n = 2k - 1$. Then the greatest length of a snake cycle in the king graph of size $n \times n$ is $(n^2 - 1)/2$. Every lift of a fewest-turn Hamiltonian cycle in the grid graph of size $k \times k$ is a longest snake cycle in the king graph of size $n \times n$. Conversely, every longest snake cycle in the king graph of size $n \times n$ can be obtained uniquely as a lift of some fewest-turn Hamiltonian cycle in the grid graph of size $k \times k$. \end{theorem}

\medskip

This time around, our analysis of paths carries over to cycles nearly verbatim.

\medskip

\begin{proof} The upper bound follows by the same argument as Theorem \ref{king-odd-path}. The lower bound is a corollary of the structure description. Finally, the structure description follows by the same argument as Theorem \ref{king-odd-structure}. \end{proof}

\medskip

With the remaining two classes, the main difficulty is this: In both Sections \ref{king-even} and \ref{king-odd-ii}, we introduce the half-sized square board $B$ of side $\lceil n/2 \rceil$ together with its grid graph $H$, and our reasoning relies heavily on the properties of the Hamiltonian paths of $H$. For Theorems \ref{king-0mod4} and \ref{king-3mod4}, it is the Hamiltonian cycles of $H$ that matter instead. When $n \equiv 1 \pmod 4$ or $n \equiv 2 \pmod 4$, however, the side of $B$ is odd and $H$ does not admit a Hamiltonian cycle.

This throws a substantial wrench in the works. While our upper bounds all go through as before, the constructions that support the lower bounds do not, and the gap which opens between the two appears to be difficult to close.

We continue with some tentative remarks.

Let $\mathcal{N}_r$ be the set of all positive integers $n$ such that $n \equiv r \pmod 4$.

Fix $n^\star$ in $\mathcal{N}_r$ with $n^\star \ge 7$. For each $n$ in $\mathcal{N}_r$ with $n \ge n^\star$, construct the subset $D_n$ of the standard board of size $n \times n$ as follows: Take all cells of the form $(x, y)$ with $x - y \ge 1$, $x + y \le n - 2$, $y$ even, and $0 \le y \le (n - n^\star)/2$ together with their images under the symmetries of the board. Delete the cells $(0, 2)$, $(0, 3)$, and $(0, 4)$, and replace them with the cells $(1, 2)$ and $(1, 4)$. Finally, for all even $i$ with $2 \le i \le (n - n^\star)/2$, delete the three cells in the set $(i, i) + \{(0, 1), (0, 3), (0, 4)\}$, and replace them with the three cells in the set $(i, i) + \{(0, 0), (1, 2), (1, 4)\}$. Note that $D_n$ is the vertex set of a snake path in $\kinggraph{n \times n}$.

We say that $\mathcal{N}_r$ \emph{crystallises} at $n^\star$ when, for all $n$ in $\mathcal{N}_r$ with $n \ge n^\star$ and every longest snake cycle $C$ in $\kinggraph{n \times n}$, there is a symmetry $\pi$ of the corresponding board such that the set of all cells of $\pi(C)$ outside of the concentric subboard of size $(n^\star - 4) \times (n^\star - 4)$ coincides with $D_n$.

Thus, in particular, if $\mathcal{N}_r$ crystallises, then there are two constants $\ell$ and $\mu$ such that the greatest length of a snake cycle in $\kinggraph{n \times n}$ is $n^2/2 - \ell$ and the number of essentially distinct snake cycles which attain this greatest length is $\mu$ for all $n$ in $\mathcal{N}_r$ with $n \ge n^\star$. Furthermore, each one of these cycles is asymmetric, and so for all such $n$ the total number of longest snake cycles in $\kinggraph{n \times n}$ is $8\mu$.

The proof of Theorem \ref{king-0mod4} shows that the class $\mathcal{N}_0$ crystallises at $n^\star = 12$ with $\ell = 1$ and $\mu = 6$.

The author finds it reasonably plausible that each one of the classes $\mathcal{N}_1$ and $\mathcal{N}_2$ might crystallise as well. Experimental data suggests that perhaps the class $\mathcal{N}_1$ crystallises at $n^\star = 13$ with $\ell = 5/2$ and $\mu = 69$ whereas the class $\mathcal{N}_2$ crystallises at $n^\star = 14$ with $\ell = 3$ and $\mu = 72$.

The next result might be helpful in the case of the class $\mathcal{N}_1$. The two cycles of Theorem~\ref{near-hamiltonian} appear to be related to the longest snake cycles of $\kinggraph{(2k - 1) \times (2k - 1)}$ in a way somewhat similar to how the fewest-turn Hamiltonian paths and cycles of grid graphs are related to the paths and cycles of Theorems \ref{king-odd-structure} and \ref{king-3mod4}.

We define a \emph{near-Hamiltonian cycle} of a graph $G$ to be a cycle in $G$ which visits all but one vertices of $G$.

\medskip

\begin{theorem} \label{near-hamiltonian} Let $k$ be an odd positive integer with $k \ge 5$. Then the smallest number of turns in a near-Hamiltonian cycle of the grid graph of size $k \times k$ is $2k$. Furthermore, for all such $k$, there are exactly $16$ near-Hamiltonian cycles which attain this smallest number. These cycles are all asymmetric, and so two of them are essentially distinct. \end{theorem}

\medskip

For completeness, there is a unique near-Hamiltonian cycle with the smallest number of turns, namely four, when $k = 3$.

Just as in Theorem \ref{king-0mod4}, the cycles of Theorem \ref{near-hamiltonian} are shaped like double spirals.

\medskip

\begin{proof} Let $B$ be the standard board of size $k \times k$, let $H$ be the grid graph on $B$, and let $C$ be a near-Hamiltonian cycle in $H$. Denote the unique cell of $B$ which $C$ omits by~$o$.

Suppose, for the sake of contradiction, that there are a row and a column of $B$ without an edge of $C$. Then the cell at their intersection must be $o$, and it cannot lie on the boundary of $B$. Since $C$ visits the neighbours of $o$ in $H$ but the row and column of $o$ do not contain edges of $C$, it follows that all edges of the cycle $(o + (1, 1)) \spring (o + (-1, 1)) \spring (o + (-1, -1)) \spring (o + (1, -1)) \spring (o + (1, 1))$ must be in $C$. However, this cycle is not near-Hamiltonian when $k \ge 5$, a contradiction.

Suppose, for concreteness, that every row contains an edge of $C$. Define the segments of $C$ as in Section \ref{king-odd-ii}. Since every row contains a horizontal segment of $C$, and the endpoints of each such segment are turns, we get that $C$ makes at least $2k$ turns altogether. Suppose, from now on, that this bound is attained and that every row contains exactly one horizontal segment of $C$. Thus, in particular, $o$ cannot lie in the lowermost or topmost row of $B$ unless it is a corner cell of $B$.

Suppose, for the sake of contradiction, that the leftmost and rightmost columns of $B$ contain one vertical segment of $C$ each. Then they cannot contain $o$ unless it is a corner cell of $B$. When, say, $o = (0, 0)$, it follows that $C$ must contain all edges of the cycle $(0, 1)$---$(1, 1)$---$(1, 0) \spring (k - 1, 0) \spring (k - 1, k - 1) \spring (0, k - 1) \spring (0, 1)$. Otherwise, when $o$ is not on the boundary of $B$, it follows that $C$ must contain all edges of the cycle $(0, 0) \spring (k - 1, 0) \spring (k - 1, k - 1) \spring (0, k - 1) \spring (0, 0)$. However, in both cases the cycle in question is not near-Hamiltonian when $k \ge 5$, a contradiction.

Suppose, for concreteness, that the leftmost column of $B$ contains at least two vertical segments of $C$. Since the number of vertical segments in $C$ is the same as its number of horizontal segments, and we have already assumed that the latter number equals $k$, we get that some column $u$ of $B$ does not contain any edges of $C$.

Consequently, $C$ crosses over $u$ every time when it visits this column. Since the total number of crossings must be even, and $u$ contains an odd number of cells, we obtain that $o$ must be in $u$.

It follows that there is exactly one column of $B$ without vertical segments of $C$. (Since each such column must contain $o$.) Thus the leftmost column of $B$ must contain exactly two vertical segments of $C$ and all other columns except for $u$ must contain exactly one vertical segment of $C$ each. (Since there are a total of $k$ vertical segments in $C$.)

Let the two vertical segments of $C$ in the leftmost column of $B$ be $(0, 0) \spring (0, w)$ and $(0, w + 1) \spring (0, k - 1)$. Suppose, for concreteness, that $1 \le w \le \lfloor k/2 \rfloor - 1$.

From this point on, we establish the identity $w = 1$ and the desired result together, by induction on $k$. The base case $k = 5$ is straightforward. For the induction step, suppose that $k \ge 7$ and that we have already settled the question on all smaller boards.

Let us delete the subpath $(1, w)$---$(0, w) \spring (0, 0) \spring (k - 1, 0) \spring (k - 1, k - 1) \spring (0, k - 1) \spring (0, w + 1)$---$(1, w + 1)$ from $C$, and let us replace it with the edge $(1, w)$---$(1, w + 1)$. The result will be a near-Hamiltonian cycle $C^\star$ in the grid graph $H^\star$ on the concentric subboard $B^\star$ of $B$ of size $(k - 2) \times (k - 2)$ given by $B^\star = [1; k - 2]^2$.

Since we have deleted at least six turns from $C$ and we have added at most two new ones in their place, $C^\star$ can make at most $2k - 4$ turns altogether. Thus our induction hypothesis applies to it, and so in fact $C^\star$ makes exactly $2k - 4$ turns, two of which are at cells $(1, w)$ and $(1, w + 1)$ where they form the subpath $(2, w)$---$(1, w)$---$(1, w + 1)$---$(2, w + 1)$. Still by the induction hypothesis, $C^\star$ omits exactly one edge of the cycle $(1, 1) \spring (k - 2, 1) \spring (k - 2, k - 2) \spring (1, k - 2) \spring (1, 1)$, and that edge is an image of the edge $(1, 2)$---$(1, 3)$ under a symmetry of $B^\star$.

We conclude that $w = 1$ and exactly two of the fewest-turn near-Hamiltonian cycles of $H^\star$ fit as a suitable $C^\star$. Therefore, there are exactly two essentially distinct fewest-turn near-Hamiltonian cycles in $H$, both of them asymmetric, and the induction step is complete. \end{proof}

%% file: snake-08-further-ii.tex
\section{Further Work on Leaper Graphs} \label{further-ii}

In this section, we collect some additional results and open problems on leaper graphs. (To be defined shortly.)

For the knight, it would be interesting to see a human-friendly proof of the upper bound in Theorem \ref{knight-all}. Or, if not that, then at least it would be nice to know if there is an unweighted knight graph with pseudosnake density $1/2$ which we could have used instead of the weighted one.

One natural direction of generalisation for our results in Section \ref{knight} is offered by leapers.

Let $p$ and $q$ be nonnegative integers with $p \le q$, not both zero. A \emph{$(p, q)$-leaper} is a fairy chess piece which moves as a generalised knight, leaping $p$ units away along one coordinate axis and $q$ units away along the other.

Let $L$ be a $(p, q)$-leaper. The \emph{leaper graph} of $L$ on a set of cells $S$, denoted $\mathcal{G}(L, S)$, is defined similarly to the king and knight graphs on $S$, except that the adjacency condition becomes $\{|x' - x''|, |y' - y''|\} = \{p, q\}$ instead.

Let $d = \gcd(p, q)$. Then the leaper graph of $L$ on the board of size $m \times n$ is the disjoint union of several isomorphic copies of the leaper graphs of a $(p/d, q/d)$-leaper on the boards of sizes $\lfloor m/d \rfloor \times \lfloor n/d \rfloor$, $\lfloor m/d \rfloor \times \lceil n/d \rceil$, $\lceil m/d \rceil \times \lfloor n/d \rfloor$, and $\lceil m/d \rceil \times \lceil n/d \rceil$, each copy scaled up by a factor of $d$. Thus we can safely assume that $d = 1$.

For $p$ and $q$ relatively prime, $L$ is known as \emph{free} when $p + q$ is odd and \emph{half-free} when it is even. Briefly, one reason for this distinction is that $\mathcal{G}(L, \mathbb{Z}^2)$ is connected when $L$ is free but consists of two connected components when $L$ is half-free.

A \emph{skew} leaper is one for which $p$ and $q$ are positive and distinct. The only non-skew leapers with relatively prime $p$ and $q$ are the $(0, 1)$-leaper, known as the \emph{wazir}, and the $(1, 1)$-leaper, known as the \emph{fers}. Of course, the wazir graph on a board coincides with the grid graph on that board. Furthermore, $\mathcal{G}(\text{Fers}, m \times n)$ can also be viewed as the direct product of two paths with $m$ and $n$ vertices, respectively.

We take a look at the wazir and fers first, and after that we will focus on skew leapers.

Even though Question \textbf{A} for grid graphs on rectangular boards is a very natural thing to ask, the only earlier reference for it known to the author as of the time of writing is \cite{R}, a puzzle game website where players are invited to construct snake paths in grid graphs on square boards, with longer paths scoring higher.

An asymptotic estimate is straightforward to obtain. The argument we give for the upper bound is not new; it is essentially identical to the argument used in \cite{MPR} to bound from above the pseudosnake density of $\gridgraph{\mathbb{Z}^2}$. (We discuss one natural way to define the pseudosnake density of certain infinite graphs below.) For the lower bound, the general strategy we outlined in Section \ref{knight} goes through without a hitch. Once again, \cite{MPR} contains the same pseudosnake in $\gridgraph{\mathbb{Z}^2}$.

\medskip

\begin{proposition} \label{wazir} Let $m$ and $n$ be positive integers. Then both the longest snake path and the longest snake cycle in the grid graph of size $m \times n$ are of length $2mn/3 + \mathcal{O}(m + n)$. \end{proposition}

\medskip

\begin{proof} Let $A$ be the standard board of size $m \times n$ and let $G$ be the grid graph on $A$.

For the upper bound, let $P$ be a snake path in $G$; the argument for cycles is similar.

Let $S$ be the vertex set of $P$ and let $T$ be the complement of $S$ within $A$. Then nearly every cell of $S$ is adjacent to two cells of $T$; each exception is either an endpoint of $P$ or near the boundary of $A$, and so there are $\mathcal{O}(m + n)$ of them. On the other hand, every cell of $T$ is adjacent to at most four cells of $S$. Thus $2|S| \le 4|T| + \mathcal{O}(m + n)$, and so $|S| \le 2|S \cup T|/3 + \mathcal{O}(m + n)$ as well.

We move on to the lower bound. Let $S_\infty$ be the set of all cells $(x, y)$ in $\mathbb{Z}^2$ with $x \not \equiv y \pmod 3$. Then $S_\infty$ induces a pseudosnake $P_\infty$ in $\gridgraph{\mathbb{Z}^2}$.

Suppose without loss of generality that $m \ge 10$ and $n \ge 10$, let $A^\star$ be the subboard of $A$ given by $A^\star = [4; n - 5] \times [4; m - 5]$, and also let $P^\star$ be the restriction of $P_\infty$ to $A^\star$. Then $P^\star$ is the disjoint union of several paths, and it is straightforward to add several cells out of $A \setminus A^\star$ to $P^\star$ so as to stitch these paths together into a single snake path or cycle. \end{proof}

\medskip

The fers can be handled similarly.

\medskip

\begin{proposition} \label{fers} Let $m$ and $n$ be positive integers. Then both the longest snake path and the longest snake cycle in the fers graph of size $m \times n$ are of length $mn/3 + \mathcal{O}(m + n)$. \end{proposition}

\medskip

The proof is analogous to that of Proposition \ref{wazir}, and we omit it.

For wazir and fers graphs on rectangular boards, it might be possible to obtain exact answers to Question \textbf{A}. By way of experimental data, \cite{R} contains a table listing the greatest length of a snake path in the grid graph of size $n \times n$ for all $n$ with $2 \le n \le 15$.

We continue on to skew leapers. Suppose, for concreteness, that $p < q$.

It seems highly likely that an exact answer to Question \textbf{A} would be out of reach for skew leapers on arbitrary rectangular boards, or even on arbitrary square boards. For this reason, we propose a weakened version of it.

\medskip

\begin{question} \label{leaper} Let $L$ be a skew leaper. What are some interesting lower and upper bounds for the greatest length of a snake path or cycle of $L$ on a given rectangular board? \end{question}

\medskip

Let us pick some of the low-hanging fruit.

When $L$ is half-free, let $\operatorname{Free}(L)$ denote the free $(p', q')$-leaper with $p' = (q - p)/2$ and $q' = (p + q)/2$.

Suppose, now, that $L$ is a skew free leaper. We will consider this case first, and then for Proposition \ref{estimate} we will reduce the half-free case to the free case using the transformation introduced above.

In all of the following asymptotic estimates, the implicit constants in the $\mathcal{O}$-terms depend on $L$.

One construction will be particularly useful to us, and so we introduce special notation for it:

Let $P$ be a pseudosnake in $\mathcal{G}(L, m \times n)$ with vertex set $S$. Consider the union $T$ of all sets of cells of the form $((n + q)i, (m + q)j) + S$, over all integers $i$ and $j$. Then the induced subgraph on vertex set $T$ is a pseudosnake in $\mathcal{G}(L, \mathbb{Z}^2)$, as the translation copies of $P$ which this subgraph consists of are too far away from one another to interact in any way. We denote this pseudosnake by $\Upsilon(m \times n, P)$.

Now let $\tau_n$ be the pseudosnake density of $\mathcal{G}(L, n \times n)$.

Observe that the sequence $\{\tau_n\}_{n = 1}^\infty$ converges. Indeed, let $n^2 \le N$. Since we can cover the board of size $N \times N$ with $\lceil N/n \rceil^2$ subboards of size $n \times n$, we get that $\tau_N \le \tau_n + \mathcal{O}(1/n)$. On the other hand, fix a largest pseudosnake $P$ in $\mathcal{G}(L, n \times n)$. Then the restriction of $\Upsilon(n \times n, P)$ to the board of size $N \times N$ is a pseudosnake in $\mathcal{G}(L, N \times N)$, and so $\tau_N \ge \tau_n + \mathcal{O}(1/n)$.

We denote $\tau(L) = \lim_{n \to \infty} \tau_n$, and we call this the pseudosnake density of $\mathcal{G}(L, \mathbb{Z}^2)$ or, for short, the pseudosnake density of $L$.

\medskip

\begin{question} \label{pseudo-density} Let $L$ be a skew free leaper. What is the pseudosnake density of $L$? Or, alternatively, what are some interesting lower and upper bounds for it? \end{question}

\medskip

Consider the four-dimensional infinite grid graph $\gridgraph{\mathbb{Z}^4}$. We define its pseudosnake density $\eta$ similarly to how we defined $\tau(L)$. Observe that $\eta$ is an absolute constant which does not depend on $L$. The pseudosnake densities of infinite grid graphs with arbitrarily many dimensions have been studied before. \cite{MPR}

\medskip

\begin{proposition} \label{density-bounds} For all skew free leapers $L$, the pseudosnake density of $L$ satisfies $1/2 \le \tau(L) \le \eta$. \end{proposition}

\medskip

\begin{proof} For the lower bound, it suffices to exhibit a doubly periodic pseudosnake in $\mathcal{G}(L, \mathbb{Z}^2)$ with density $1/2$.

Since $L$ is free, exactly one of $p$ and $q$ is even. Denote that even value by $r$. When $r \equiv 2 \pmod 4$, let $S$ be the set of all cells $(x, y)$ such that $\lfloor x/2 \rfloor$ is even. Otherwise, when $r \equiv 0 \pmod 4$, let $S$ be the set of all cells $(x, y)$ such that $\lfloor x/2 \rfloor + y$ is even. Then the induced subgraph of $\mathcal{G}(L, \mathbb{Z}^2)$ on vertex set $S$ works.

For the upper bound, let $n$ be a positive integer, fix a largest pseudosnake $P$ in $\mathcal{G}(L, n \times n)$, and let $Q = \Upsilon(n \times n, P)$. Then $Q$ is a doubly periodic pseudosnake in $\mathcal{G}(L, \mathbb{Z}^2)$ with density $\tau_n + \mathcal{O}(1/n)$.

Consider, now, the induced subgraph $R$ of $\gridgraph{\mathbb{Z}^4}$ whose vertex set consists of all four-dimensional integer points $(x_1, x_2, x_3, x_4)$ such that $x_1(p, q) + x_2(q, p) + x_3(-p, q) + x_4(-q, p)$ is a cell of $Q$. Then $R$ is a quadruply periodic pseudosnake in $\gridgraph{\mathbb{Z}^4}$ with the same density as $Q$. We let $n$ grow without bound, and the conclusion follows. \end{proof}

\medskip

Had it been the case that $\eta = 1/2$, Proposition \ref{density-bounds} would have resolved Question~\ref{pseudo-density} immediately. However, it has been demonstrated that $649/1296 \le \eta \le 20/39$, with $649/1296 \approx 0.50077$ and $20/39 \approx 0.51282$. \cite{MPR} Still, Proposition \ref{density-bounds} and this result together imply, for all skew free leapers $L$, that $1/2 \le \tau(L) \le 20/39$. Since the gap between these two bounds is rather narrow, and the graph $\mathcal{G}(L, \mathbb{Z}^2)$ is, in some intuitive sense, more crowded than $\gridgraph{\mathbb{Z}^4}$, it seems plausible that in fact $\tau(L) = 1/2$ for all skew free leapers $L$. As we saw in Section \ref{knight}, this is certainly true of the knight.

The natural connection between pseudosnake density and snake paths and cycles is as follows:

\medskip

\begin{proposition} \label{estimate} Let $L$ be a skew leaper and let $m$ and $n$ be positive integers. When $L$ is free, the greatest length of a snake path or cycle of $L$ on the board of size $m \times n$ does not exceed $\tau(L) \cdot mn + \mathcal{O}(m + n)$. Furthermore, when $L$ is half-free, it does not exceed $\tau(L)/2 \cdot mn + \mathcal{O}(m + n)$. \end{proposition}

\medskip

\begin{proof} Let $P$ be a snake path or cycle in $\mathcal{G}(L, m \times n)$ with $s$ cells.

Suppose first that $L$ is free. Since $\Upsilon(m \times n, P)$ is a doubly periodic pseudosnake in $\mathcal{G}(L, \mathbb{Z}^2)$, its density $s/(m + q)(n + q)$ does not exceed $\tau(L)$.

Suppose, now, that $L$ is half-free. Then $\varepsilon = (x + y) \bmod 2$ is constant over all cells $(x, y)$ of $P$. Instead of $\Upsilon(m \times n, P)$, take its image under the transformation $(x, y) \to ((x - y + \varepsilon)/2, (x + y + \varepsilon)/2)$. This is a pseudosnake in $\mathcal{G}(\operatorname{Free}(L), \mathbb{Z}^2)$, and from this point on the argument continues as before. \end{proof}

\medskip

The author finds it reasonably plausible that the upper bounds of Proposition \ref{estimate} might in fact be attained for all skew leapers. As we saw in Section \ref{knight}, this is indeed the case for the knight.

Note that our construction for the lower bound in the proof of Proposition \ref{density-bounds} yields a doubly periodic pseudosnake in $\mathcal{G}(L, \mathbb{Z}^2)$ where all cells are of degree exactly two and there are no finite cycles. Thus for free leapers $L$ with $\tau(L) = 1/2$ and their corresponding half-free leapers, this construction might play a role in a proof that the upper bounds of Proposition \ref{estimate} are attained which follows some variant of the strategy we outlined in Section \ref{knight}.

%% file: snake-09-ack.tex
\section*{Acknowledgements} \label{ack}

The author would like to thank Professor Donald Knuth for introducing him to the subject of the longest snake paths and cycles in chess piece graphs.

%% file: snake-11-add.tex
\section*{Addendum} \label{add}

Since the completion of this paper, the author has become aware of additional relevant references.

Kolmogorov Cup, 2008, round 2, major league, problem 2, by Alexandre Chapovalov and Maxim Chapovalov, \mbox{\url{https://turmath.ru/kolm/files/archive/kolm12.zip}}, is about the greatest length of a snake path in $\kinggraph{100 \times 100}$. (First league, problem 2, second league, problem 1, and major junior league, problem 5 are versions with $\kinggraph{8 \times 8}$.)

Mathematical Festival, 2002, grade 7, problem 5, by Igor Akulich, \url{https://olympiads.mccme.ru/matprazdnik/image/matpr2002.pdf}, together with its official solution is about the greatest length of a snake path in $\gridgraph{8 \times 8}$. (Grade 6, problem 4 is a version with $\gridgraph{7 \times 7}$.)

Two formulas of OEIS entry A331968, added in 2020 by Elijah Beregovsky, \url{https://oeis.org/A331968}, with proofs in the entry's revision history, bound the greatest length of a snake path in $\gridgraph{n \times n}$ within an interval of size $\mathcal{O}(n)$.